\documentclass[11pt]{amsart}
\usepackage{pdfsync}
\usepackage{a4,fullpage}
\usepackage{enumerate}
\usepackage{graphicx}
\usepackage{latexsym}
\usepackage{amsfonts}
\usepackage[utf8]{inputenc}
\usepackage{amsthm}
\usepackage{amsmath}
\usepackage{amssymb}
\usepackage[
            unicode,
            breaklinks=true,
            colorlinks=true]{hyperref}
\usepackage[active]{srcltx}
\usepackage{amscd}
\usepackage[all]{xy} 
\usepackage{mathrsfs} 
\usepackage{stmaryrd} 
\usepackage{amsopn} 
\usepackage{eepic}
\usepackage{epic}
\usepackage{eucal}
\usepackage{epsfig}
\usepackage{psfrag}
\usepackage{amscd}
\usepackage{mathrsfs}
\usepackage{stmaryrd}
\usepackage{amsopn}
\newtheorem{theorem}{Theorem}[section]

\newtheorem{problem}[theorem]{Problem}
\newtheorem{coreproblem}[theorem]{Core Problem}
\newtheorem{example}[theorem]{Example}
\newtheorem{question}[theorem]{Question}
\newtheorem{classproblem}[theorem]{Classification Problem}

\newtheorem{conjecture}[theorem]{Conjecture}

\newtheorem{definition}[theorem]{Definition}
\newtheorem{remark}[theorem]{Remark}

\newcommand{\Hilbert}{{\mathcal{H}}}  
\newcommand{\RM}{{\mathbb R}}
\newcommand{\R}{{\mathbb R}}
\renewcommand{\C}{{\mathbb C}}
\newcommand{\Z}{{\mathbb Z}}

\newcommand{\N}{{\mathbb N}}
\newcommand{\NM}{{\mathbb N}}

\newcommand{\Cinf}{C^\infty}
\newcommand{\got}[1]{\mathfrak{#1}}

\newcounter{todo}
\newcommand{\op}[1]{\!\!\mathop{\rm ~#1}\nolimits}

\newcommand{\scriptop}[1]{\!\!\mathop{\mbox{\rm \scriptsize ~#1}}\nolimits}

\newcommand{\ii}{\mathrm{i}}

\newcommand{\om}{\omega}

\newcommand{\De}{\Delta}

\renewcommand{\leq}{\leqslant}

\begin{document}

\title[Symplectic and spectral theory of integrable systems]{First steps in symplectic and spectral theory of
    integrable systems}

\author{\'Alvaro
  Pelayo\,\,\,\,\,\,\,\,\,\,\,\,\,\,\,\,\,\,\,\,\,\,\,\,\,\,San V\~{u}
  Ng\d{o}c}

\maketitle

\begin{abstract}
  The paper intends to lay out the first steps towards constructing a
  unified framework to understand the symplectic and spectral theory
  of finite dimensional integrable Hamiltonian systems.
    While it is difficult to know what the best approach to such a large
  classification task would be, it is possible to single out some
  promising directions and preliminary problems. This paper discusses
  them and hints at a possible path, still loosely defined, to arrive
  at a classification. It mainly relies on recent progress concerning
  integrable systems with only non-hyperbolic and non\--degenerate
  singularities.

  This work originated in an attempt to develop a theory aimed at
  answering some questions in quantum spectroscopy. Even though
  quantum integrable systems date back to the early days of quantum
  mechanics, such as the work of Bohr, Sommerfeld and Einstein, the
  theory did not blossom at the time.  The development of
  semiclassical analysis with microlocal techniques in the last forty
  years now permits a constant interplay between spectral theory and
  symplectic geometry. A main goal of this paper is to emphasize the
  symplectic issues that are relevant to quantum mechanical integrable
  systems, and to propose a strategy to solve them.
\end{abstract}

\section{Program outline} \label{sec:intro}

This paper suggests an approach to work towards a \emph{symplectic and
  spectral} classification of finite dimensional integrable
Hamiltonian systems.

It is a bottom to top path, consisting of many problems, in which each
one tries to deal with only one major difficulty. Some of these
problems are labelled as ``Core Problems'' to indicate that they are
essential milestones towards the so-called ``Classification Problems''
that are the main incentives of this paper.  Our proposed strategy to
deal with each of these problems is inspired by our articles
\cite{PeVN2009,PeVN2011}. In order to facilitate the reading of the
paper we briefly review the main results and proofs of these articles,
as well as some essential ingredients from the literature used in
them.

We are optimistic about the potential of the ideas we present in this
paper, while we recognize many technical and conceptual challenges to
implement them. Of course, we do not know if some of the challenges
could be impossible to surpass.  Two key tools in the modern
approaches to integrable systems are:
\begin{itemize}
\item[(i)] the study of \emph{singular affine structures}, and
\item[(ii)] the \emph{symplectic linearization theorems} for
  non\--degenerate singularities\footnote{In the
    $\op{C}^{\infty}$\--smooth category, we expect major difficulties
    in the cases where the integrable systems have degenerate
    singularities (this is not so in the algebraic setting, as seen by
    the works of Garay and van Straten).  The lack of symplectic
    linearization theorems for degenerate singularities makes their
    study challenging.  At this time we neither have much information
    about such degenerate integrable systems, nor do we have methods
    to analyze them (see Section \ref{sec:degenerate} for further
    comments). Algebraic geometry seems a natural setting for the
    study degenerate singularities.  The case where a system has
    hyperbolic singularities is rich and interesting, as we know from
    examples by Zhilinski\'{\i} and a few results by V\~{u} Ng\d{o}c,
    Dullin, Zung, and Bolsinov, among others. In addition, a key tool
    on the spectral theory side is the microlocal analysis of Toeplitz
    operators developed in the past decade.}.
\end{itemize}
These tools dominate, implicitly or explicitly, the proof strategies
that we outline.

In this framework, the known integrable systems should be understood
in terms of a collection of invariants.  Six motivational examples of
the program are: the \emph{coupled spin\--oscillator}, the
\emph{spherical pendulum}, the \emph{Lagrange top}, the
\emph{two\--body problem,} the \emph{Kowalevski top}, and the
\emph{three wave interaction}; these examples are explained in more
detail in Section \ref{sec:first}.

Some of these examples have hyperbolic or degenerate singularities.
Analyzing the case of integrable systems which have hyperbolic
singularities, will probably involve major breakthroughs (see Section
\ref{sec:hyperbolic}); it remains largely unexplored.  Our proposal
for a unifying approach aims at identifying the essential features of
integrable systems through the computation of invariants that classify
them. In this sense, a success in our approach would help to reconcile
the vast amount of literature for specific systems, with the
theoretical approaches developed in recent times (see Section
\ref{sec:comp}).

The present paper \emph{does not} intend to be a survey, but rather a
fast description of a few preliminary ideas which fit into a larger
plan. In this sense it is more a ``work guide'' than a research paper.
It is evident from the very few articles we cite in the paper, that
our choice of references is mostly a practical one. We refer to
Section \ref{sec:r} for further discussion in this direction, and for
references to survey papers where more extensive bibliographies may be
found.

The paper could be read independently, but it is probably best read
simultaneously or after our article \cite{PeVN2012} which gives a
succinct overview of the current research in the subject, with
references to previous works. Both papers complement each other and
have different focuses.

\subsection{\textcolor{black}{Origins}}

The first steps in the approach to classifying integrable systems
advocated in this paper originated in our attempts to develop
technology to answer inverse type questions about quantum systems in
molecular spectroscopy.

While the notion of a quantum integrable system dates back to the
early times of quantum mechanics (e.g., in the works of Bohr,
Sommerfeld, and Einstein), the basic results in the symplectic theory
of classical integrable systems could not have been used in
Schr\"odinger's quantum setting for the simple reason that the analysis
of pseudodifferential operators in phase space was developed later. In
addition to pseudodifferential operators, this paper advocates the use
of microlocal analysis of Toeplitz operators to analyze quantization
and inverse spectral problems on arbitrary phase spaces (not
necessarily cotangent bundles).

In \cite{PeVN2009,PeVN2011}, we carried out the symplectic part of
this ``program'' for non\--hyperbolic systems with two degrees of
freedom in four dimensional phase space for which one component of the
system is periodic and proper -- these systems are called
\emph{semitoric systems}. A physical example of semitoric system is
the \emph{coupled spin\--oscillator}, mentioned above.

\subsection{\textcolor{black}{Goals}} \label{go}

\subsubsection{\textcolor{black}{Symplectic Geometry}}

Recall that a classical \emph{integrable system} is given by the
following data: (1) a $2n$\--dimensional smooth manifold $M$ equipped
with a symplectic form $\omega$, and (2) $n$ smooth
functions $$f_1,\ldots,f_n \colon M \to \mathbb{R}$$ which generate
Hamiltonian vector fields $\mathcal{H}_{f_1},\, \ldots,\,
\mathcal{H}_{f_n}$ (where $\mathcal{H}_{f_i}$ is defined by
\emph{Hamilton's equation} 
$\omega (\mathcal{H}_{f_i},\, \cdot)={\rm  d}f_i$ 
for every $1 \leq i \leq n$) that are linearly independent at
almost every point of $M$ and which pairwise commute in the sense that
the Poisson brackets vanish:
$$
\{f_i, \, f_j\}:=\omega(\mathcal{H}_{f_i},\,\mathcal{H}_{f_j})=0,
\quad \textup{for all} \quad 1 \leq i,\,j \leq n.
$$

\begin{figure}[h]
  \centering
  \includegraphics[width=0.43\linewidth]{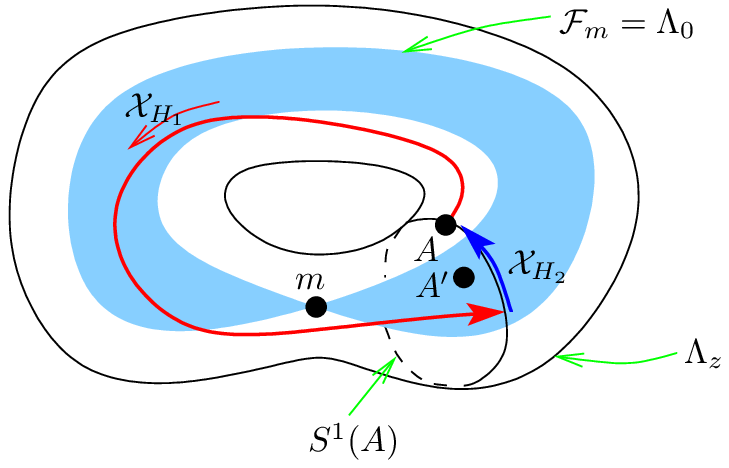}
  \caption{The most interesting features of integrable systems are in
    their singularities. The singularities encode both local and
    global information of the dynamics, geometry, and topology of the
    system.  The picture shows a singular foliation near a singular
    leaf $\Lambda_0$ of an integrable system $(H_1,\,H_2) \colon
    \mathbb{R}^4 \to \R^2$ given by two Hamiltonians $H_1 \colon
    \mathbb{R}^4 \to \mathbb{R}$ and $H_2 \colon \mathbb{R}^4 \to
    \mathbb{R}$. The point $m$ is a singularity of so called
    \emph{focus\--focus} type. The point $A$ denotes a regular point
    close to $m$, which lies in a regular fiber $\Lambda_c$ of the
    system $(H_1,\,H_2) \colon \mathbb{R}^4 \to \mathbb{R}$.  The
    Hamiltonian $H_1$ generates a vector field $\mathcal{H}_{H_1}$
    with a periodic flow, while the other Hamiltonian generates a
    vector field $\mathcal{H}_{H_2}$ with a hyperbolic and
    non\--periodic flow. The set $S^1(A)$ denotes a periodic orbit of
    $H_1$.}
  \label{fig:pp}
\end{figure}

We also refer to the map $F:=(f_1,\, \ldots,\,f_n) \colon M \to
\mathbb{R}^n$ as the \emph{integrable system} itself.  Although the
path to arrive at it is still loosely defined, the first goal of the
program is clear.
\\
\\
{\bf \emph{Symplectic Goal of Program for Integrable Systems}}:
\emph{Prove that large classes of finite dimensional integrable
  Hamiltonian systems $ F=(f_1,\, \ldots,\, \, f_n) \colon M \to
  \mathbb{R} $ are characterized, up to symplectic isomorphisms, by an
  explicit list of invariants} (a program proposal to arrive at this
result is outlined in section \ref{sec:symplectic}).

\subsubsection{\textcolor{black}{Spectral Geometry}}

The driving force behind this paper, and behind the goal above, is the
inverse question: if we know the spectrum of a quantum physical system
which is integrable (see Figure \ref{fig:spectra2}), can we
reconstruct the classical system from it?  Answering this question
involves studying symplectic and spectral invariants. The study of
these invariants is the common theme of the present paper. It is
connected, at a fundamental level, with the analysis of the
singularities of integrable systems which encode information about the
global behavior of the system (Figure \ref{fig:pp}).

>From a mathematical point of view, this question leads to the second
goal of this program: to develop an inverse spectral result along the
following lines. A \emph{quantum integrable system} is defined as a
family of commuting operators $T^{\hbar}_1,\,\ldots, \, T^{\hbar}_n$
on Hilbert spaces $\mathcal{H}_{\hbar}$, indexed by $\hbar \rightarrow
0$, whose principal symbols $f_1,\, \ldots,\, f_n$ form a classical
integrable system on a symplectic manifold $(M,\, \omega)$. The
\emph{semiclassical joint spectrum} of the system is given by the
collection of joint spectra of $T^{\hbar}_1,\, \ldots,\, T^{\hbar}_n$
(see Figure \ref{fig:spectra2} for one such element.)  Next us explain
the terminology more concretely when $M$ is a \emph{compact manifold}
(\emph{but} the spectral goal refers to all manifolds, compact or
not).

In the \emph{compact} case the commuting operators which we consider
are ``Toeplitz operators", and the quantization of the symplectic
manifold is the ``geometric quantization". Note that the
\emph{Spectral Goal} assumes that we have a quantization of the
\emph{manifold}, and the existence of such a quantization is far from
obvious. A now standard procedure by B. Kostant and J.M. Souriau in
the 1960s is to introduce a \emph{prequantum} bundle $\mathcal{L}
\rightarrow M$, that is a Hermitian line bundle with curvature
$\frac{1}{i} \om$ and a complex structure $j$ compatible with
$\omega$. One then defines the quantum space as the
space $$\mathcal{H}_k:=\mathrm{H}^0(M,\mathcal{L}^k)$$ of holomorphic
sections of tensor powers $\mathcal{L}^k$ of $\mathcal{L}$. Here the
semiclassical parameter is $\hbar=1/k$.  The parameter $k$ is a
positive integer, the semi-classical limit corresponds to the large
$k$ limit.  Associated to such a quantization there is an algebra
$\mathscr{T} (M,{\mathcal{L}} , j)$ of operators, called Toeplitz
operators.  A \emph{Toeplitz operator} is any sequence $(T_k \colon
\Hilbert_k \rightarrow \Hilbert_k )_{k \in \mathbb{N}^*}$ of operators
of the form
\begin{eqnarray} \nonumber \Big(T_k = \Pi_k f(\cdot,\, k) + R_k
  \Big)_{k\in \mathbb{N}^*} \label{qq}
\end{eqnarray}
where $f(\cdot,\,k)$, viewed as a multiplication operator, is a
sequence in $\Cinf( M)$ with an asymptotic expansion $f_0 + k^{-1} f_1
+ \ldots $ for the $\Cinf$ topology, and the norm of $R_k$ is
$\mathcal{O}(k^{-\infty})$. The algebra of Toeplitz operators plays
the same role as the algebra of semiclassical pseudodifferential
operators for a cotangent phase space.  A Toeplitz operator has a
\emph{principal symbol}, which is a smooth function on the phase space
$M$. If $T$ and $S$ are Toeplitz operators, then
$(T_k+k^{-1}S_k)_{k\in\NM^*}$ is a Toeplitz operator with the same
principal symbol as $T$.  If $T_k$ is Hermitian (i.e. self-adjoint)
for $k$ sufficiently large, then the principal symbol of $T$ is
real-valued.  Two Toeplitz operators $(T_k)_{k \in \N^*}$ and
$(S_k)_{k \in \N^*}$ \emph{commute} if $T_k$ and $S_k$ commute for
every $k$.
\\
\\
{\bf \emph{Spectral Goal of Program for Integrable Systems}}:
\emph{Show that the semiclassical joint spectrum of a quantum
  integrable system recovers the classical integrable system, up to
  symplectic isomorphisms} (a program proposal to arrive at this
result is outlined in Section \ref{sec:spectral}).
\\
\\
As we will see later, to complete this second goal, partial success
with the first goal is necessary. Put differently, achieving the
second goal for quantum integrable systems (of a certain type) passes
through the solution of the first goal for classical integrable
systems (of the type corresponding to the quantum system at hand).

\begin{figure}[h!]
  \centering
  \includegraphics[width=0.8\textwidth]{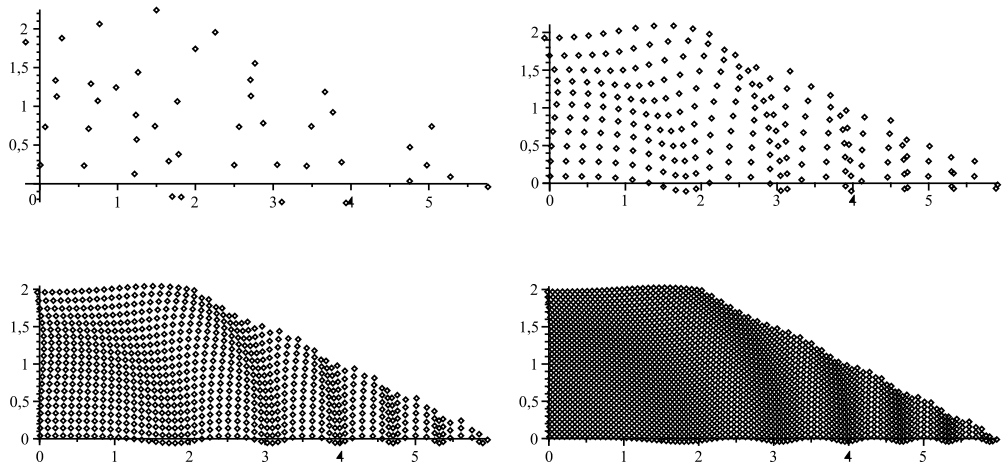}
  \caption{Sequence of images of the spectra of a quantum toric
    integrable systems as $\hbar$ goes to $0$.}
\end{figure}

\subsection{\textcolor{black}{Approach comparison}}
\label{sec:comp}

We put forward some ideas to unify the existing knowledge of concrete
integrable systems. The final goals were presented in Subsection
\ref{go}. Our approach differs from the traditional one, where
concrete individual or families of integrable systems are studied. We
do not intend to replace this classical method but to complement it by
drawing attention to an additional research focus.  For instance, the
coupled spin\--oscillator system fits in the classification theory
developed in our papers \cite{PeVN2009,PeVN2011}. However, these
articles do \emph{not} determine the invariants of the coupled
spin\--oscillator, which is a difficult computational problem
addressed in our later paper \cite{PeVN2010}; all invariants of the
coupled spin\--oscillator are found, with the exception of one for
which only a rigorous linear approximation is given (this invariant is
a formal power series).

The existing literature contains sophisticated techniques for
computing features of integrable systems explicitly.  These methods,
and their conclusions, are spread over many thousands of papers.  Our
theory for semitoric systems (developed in \cite{PeVN2009,PeVN2011})
and its extension proposed here, neither improve these traditional
approaches nor replace the invariants they compute.  Although
superficially, many of these characteristic traits seem unrelated to
our symplectic invariants, they are, of course, linked to
them. Whenever a property of an integrable system does not change
under symplectic isomorphisms which preserve the system, it must be
encoded in the finite list of invariants we construct.  Our initial
investigations have relied on existing literature and have been
motivated by it.  The ideas suggested in this paper go hand in hand
with the traditional work and are not independent from it.

\subsection{\textcolor{black}{Strategies and proofs}}

The purpose of this paper is to present many open problems, not to
outline avenues for their solution. The reason for this is
two-fold. First, we have not solved these problems ourselves and it is
very likely that surprises may lurk in the technical details. Second,
we believe that at least some of them may have proofs that would
follow the \emph{outline} of those in \cite{PeVN2009, PeVN2011}:
construct local symplectic invariants (analogues of the invariants
defined in these papers) and then symplectically glue these
pieces. The proofs of our classification theorems in \cite{PeVN2009,
  PeVN2011} -- for the so called \emph{semitoric systems} -- are
divided into separate steps, so it is easy to follow and understand
the overall approach. We review in Subsection \ref{sec:semitoric} the
strategy to prove these classification theorems (Theorem
\ref{inventiones} and Theorem \ref{acta}), emphasizing which methods
therein could be applicable in more general cases. In addition, we
discuss the difficulties encountered in each part of the overall
setup. The arguments proper to each step in \cite{PeVN2009, PeVN2011}
depend on the symplectic invariants given in these articles.

The building of these invariants is somewhat ``rigid" and is designed
for the type of systems classified in \cite{PeVN2009, PeVN2011} (i.e.,
semitoric systems).  Slightly changing the assumption on the type of
integrable system can invalidate major parts of the original
construction of the invariant (for semitoric systems).  How the
assumptions of these theorems (i.e., the type of integrable system
under consideration) can be weakened corresponds to the titles of the
subsections \ref{sec:hyperbolic}, \ref{sec:proper}, \ref{sec:higher},
\ref{sec:antiperiodic}, \ref{sec:degenerate}, and \ref{sec:topology}:
\emph{hyperbolic systems, non-proper systems, higher dimensional
  systems, non-periodic systems, degenerate systems, and independent
  topological and geometric questions} (see Section \ref{sec:first}
for examples which fit into one of these categories, or in the overlap
of several categories.)

The goal of subsections \ref{sec:hyperbolic}, \ref{sec:proper},
\ref{sec:higher}, \ref{sec:antiperiodic}, \ref{sec:degenerate}, and
\ref{sec:topology} is to outline, and briefly explain, the major
difficulties we foresee, based on our previous experience dealing with
integrable systems.  Then, in subsection \ref{sec:collaborative} (with
the title ``Collaborative efforts"), we show how to gather all the
information gleaned from the previous subsections \ref{sec:proper},
\ref{sec:higher}, \ref{sec:antiperiodic}, \ref{sec:degenerate}, and
\ref{sec:topology} to achieve a general classification.  In our
opinion, each of these subsections presents a challenge on its own
because they address somewhat unrelated technical difficulties; so we
suggest that they be considered individually.  Precisely because of
the likely ``independence" of these assumptions, putting all of them
together in Section \ref{sec:collaborative} should not be a major
technical problem. Instead, we expect it to be a book-keeping problem:
if all goes well, the general answer in subsection
\ref{sec:collaborative} would be a ``direct sum'' of the answers
provided in subsections \ref{sec:hyperbolic}, \ref{sec:proper},
\ref{sec:higher}, \ref{sec:antiperiodic}, \ref{sec:degenerate}, and
\ref{sec:topology}.

\subsection{\textcolor{black}{Directions}}

At the time of the writing of this paper, most of the problems we
state are open. We have begun work on a few of them, and where this is
the case, we explicitly say so.

The widely open directions outlined in this paper concern the
symplectic theory of systems with degenerate singularities, systems
with hyperbolic singularities, and the corresponding spectral theory
for these systems; this is explained in subsections
\ref{sec:hyperbolic}, \ref{sec:degenerate}, and parts of
\ref{sec:spectral}. There are also many open questions of geometric
and topological nature that fit within the program; some of these can
be found in subsection \ref{sec:topology}. The nonperiodic case,
treated in subsection \ref{sec:antiperiodic}, is also mostly open and
our initial investigations indicate the presence of challenges already
in dimension four. Subsection \ref{sec:collaborative} is open, but it can be undertaken
only after all the problems of the previous subsections have been
solved.

\subsection{\textcolor{black}{Timeline}}

We are writing this paper with the intention of raising awareness of
what we believe is a known fact among some specialists on integrable
systems: the time is ripe for great advances. The reason is that
several new effective tools and ideas are now understood, which was
not the case, say, ten years ago.  As already mentioned, there are
three key ingredients that are understood today much better: singular
affine structures, symplectic singularity theory, and the microlocal
analysis for Toeplitz operators. The technology and ideas developed by
many since the 1970 and 1980s, but particularly in the past ten years,
could be greatly developed with a view towards a better understanding
of the symplectic and spectral aspects of integrable systems.

We are optimistic about the chances of success of the ideas presented
here, judging from our own investigations in the past few years.
Nevertheless, we admit that there are still many technical and
conceptual challenges that need resolution and do not know if some of
our scenarios could lead to a dead-end. The program is at the
beginning stages, we explain its success highlights so far, and
outline expected challenges.  It is difficult to predict how long it
will be till major progress has been made.

We hope that this paper may be of help to readers who wish to do
research on some part of the program.

\subsection{\textcolor{black}{How to read this paper}}

The paper could be read independently, but it is probably better to
consult simultaneously our review \cite{PeVN2012}, or even to
familiarize oneself first with it. In \cite{PeVN2012} we gave an
overview of the state of the art in the subject. Both articles
together could serve as a speedy way to be immersed in the subject and
have easy access to some prominent problems.  Also, the paper need not
be read linearly.

Section \ref{sec:symplectic} deals with the symplectic theory of
integrable systems and proposes a self\--contained program on its
own. Section \ref{sec:spectral} addresses the spectral theory of
integrable systems, but it relies heavily on Section
\ref{sec:symplectic}.

The plan suggested in the present paper is dynamic and we expect that
the approaches to some of the problems, and the strategies we
describe, will change as the ideas presented here evolve, and new ones
arise.  In this sense, this paper is far from giving a definite
approach to the problems that it proposes.

\section{Symplectic theory of integrable
    systems} \label{sec:symplectic}

The unifying topic of this section is \emph{symplectic geometry},
which is a mathematical language that clearly and concisely formulates
many problems in theoretical physics. It also provides a framework in
which classical and quantum physics are treated simultaneously.
Symplectic geometry is closely connected with many areas of
mathematics. Within symplectic geometry, we focus on integrable
systems which are a fundamental class of ``explicitly" solvable
dynamical systems of interest in classical dynamics, semiclassical
analysis, partial differential equations, low-dimensional topology,
algebraic geometry, representation theory, and theoretical physics.

Many integrable systems are found in simple physical models of
classical and quantum physics. Two examples are the spherical pendulum
(see Figure \ref{fig:pendulum}) and the quantum coupled spin
oscillator (see Figure \ref{fig:spectra2}).  Integrable systems also
appear in other applied fields, such as locomotion generation in
robotics, elasticity theory, plasma physics, and planetary mission
design.  In spite of being, in some sense, ``solvable", they exhibit a
rich behavior. For instance, the symplectic dynamics around the so
called focus-focus singularities is highly complex (it has infinitely
many symplectic invariants, for example).

\subsection{\textcolor{black}{Results, tools, and
    methodology}} \label{sec:first}

\begin{figure}[h]
  \begin{center}
    \includegraphics[width=0.7\textwidth]{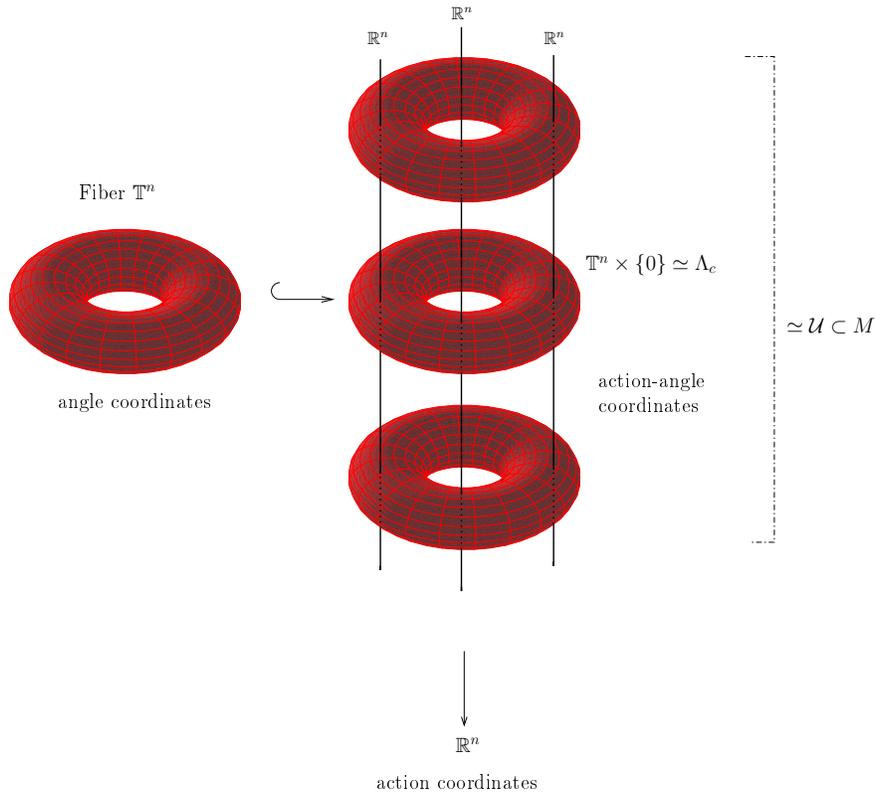}
    \caption{According to the Arnold\--Liouville\--Mineur theorem, a
      tubular neighborhood $\mathcal{U}$ of a regular fiber embeds
      symplectically into $\op{T}^*(\mathbb{T}^n)\simeq \mathbb{T}^n
      \times \mathbb{R}^n$, and the integrable system on this
      neighborhood is given by the $n$ canonical projections to
      $\mathbb{R}$.  Hence the dynamics around a regular fiber is
      simple.}
    \label{figure1}
  \end{center}
\end{figure}

Our initial motivation for studying integrable systems was to develop
a theory capable of answering questions raised by physicists and
chemists, working on quantum molecular spectroscopy. The main question
is simple to state: can a classical integrable physical system be
reconstructed from the spectrum of its (semi-classical) quantization?
Mathematically, this is a fascinating problem. The answer seems to be
``yes" for a number of cases. For instance, it is the case for toric
systems. Not only symplectic geometry is the tool to answer this
question, but it also serves as a method to help understand, even
predict, features of the spectral theory of systems in quantum
chemistry and quantum physics.

Studying this inverse question in particular, and spectral geometry of
quantum integrable systems in general, requires introducing and
combining a number of mathematical tools; this is what this project is
all about.
\\
\\
{\bf \emph{A. Previous Results.}} The behavior of many complicated
integrable systems is well understood. However, few \textit{general}
results are known. Among these, the Arnold-Liouville action angle
theorem stands out as a classical result to be found in any text on
mechanics.  In addition, there have been some remarkable successes
describing local and semilocal behavior of integrable systems such as
Eliasson's linearization theorem for non-degenerate singularities or
the work of V\~{u} Ng\d{o}c and Zung on singularities. Shortly before
Eliasson proved his theorems in the mid 1980s, Duistermaat formulated
his global action-angle obstruction theorem which remains, to this
day, one of the landmark results in the theory of integrable systems.

>From a topological point of view, the work of the Fomenko school on
the classification of of singular Lagrangian fibrations is ground
breaking. These results completely classify the topology of
non\--degenerate two degrees of freedom systems, give a method for
handling finite dimensional integrable systems in any dimension, and
sometimes even include global descriptions. The important results on
compact Hamiltonian group actions by Kostant, Atiyah, Guillemin,
Sternberg, Kirwan, and Delzant, provide another model and motivation
for the program proposed in this paper.
\\
\\
{\bf \emph{ B. New Tools.} } A number of mathematical developments to
which many people have contributed, principally since the 1970/1980s,
but particularly in the past fifteen years, have provided effective
tools that can be used to study global aspects of integrable systems.
These include methods from functional analysis and partial
differential equations (particularly, geometric techniques involving
Fourier integral operators, Toeplitz operators, linearization
theorems), symplectic geometry (developments on Lie theory, monodromy,
symplectic actions, singular reduction, relations between quantization
and reduction, global action-angle coordinates), and algebraic
geometry (singular Lagrangian fibrations, affine structures, and the
development of the theory of toric varieties).
\\
\\
{\bf \emph{C. Program methodology.}} This program relies on the
symplectic local and semilocal theory of integrable systems as well as
the detailed description of several integrable systems available in
the mathematics and physics literature.  We propose to investigate
global patterns using the aforementioned tools.  With their aid, one
can introduce a framework to describe complicated phenomena in four
dimensions, i.e., Hamiltonian integrable systems with two degrees of freedom, for
which one component of the system is periodic; these are called
\emph{semitoric integrable systems}. This program outlines a plan for
developing these results, with the ultimate goal of achieving a full
understanding of the symplectic and spectral theory of finite
dimensional integrable systems.
\\
\\
{\bf \emph{D. Motivational examples.}} The following are famous
examples of finite dimensional integrable systems which fit in the
scope of our program. The reader may revisit them after reading each
of the upcoming sections of the paper. They have been thoroughly
studied from many points of view.  Examples (1), (2), (3), (4), and
(5) are integrable systems with two degrees of freedom on
four\--dimensional phase space. Example (6) is an example of an
integrable system with three degrees of freedom on a six\--dimensional
phase space. All of the examples below are on a non\--compact phase
space.

All sections of our program, with the exception of subsection
\ref{sec:degenerate} (on degenerate systems), use the linearization
theorem of Eliasson for singularities. Hence one has to check the
non\--degeneracy of the singularities (see Definition
\ref{def:non-deg}).  We have verified this for some of the examples
below.  One could probably deduce this non\--degeneracy (or lack
thereof) from the existing literature, but we are not aware of a
reference where this has been explicitly done.

\begin{itemize}
\item[(1)] \emph{Coupled spin\--oscillator}. This integrable system
  fits into the developed and finalized theory in subsection
  \ref{sec:semitoric} (semitoric systems). We checked this in
  \cite[Section 2]{PeVN2010}. In this article there is also a detailed
  study of the coupled spin\--oscillator and its quantum counterpart.
  This example is of primary importance in physics, where it is called
  the \emph{Jaynes\--Cummings model}. There are natural extensions of
  this model to arbitrary dimensions.

\begin{figure}[h]
  \centering
  \includegraphics[width=0.25\textwidth]{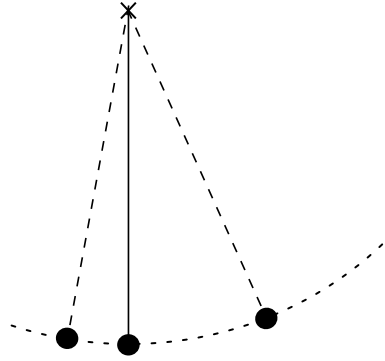}
  \caption{The spherical pendulum, a simple example of a semitoric
    system.}
  \label{fig:pendulum}
\end{figure}

\item[(2)] \emph{Spherical pendulum}.  This integrable system fits
  into subsection \ref{sec:proper} (non\--proper systems).  We checked
  this in \cite[Section 5]{PeRaVN2011}. (There is a recent paper by
  Dullin that computes invariant (ii) in Theorem \ref{inventiones} for
  the spherical pendulum).

\item[(3)] \emph{Two\--body problem.}  This integrable system fits
  into a combination of subsection \ref{sec:hyperbolic} (hyperbolic
  systems) and subsection \ref{sec:proper} (non\--proper systems).
  One of its components generates a periodic flow, but it is not a
  proper map. The singular Lagrangian fibration defined by the system
  is not proper either (because the bifurcation set and the critical
  set do not coincide).  This integrable system has hyperbolic
  singularities.

  We have not checked whether all the singularities of this system are
  non\--degenerate.  If there are degenerate singularities, then the
  study of this example also would overlap with subsection
  \ref{sec:degenerate} (degenerate systems), in addition to subsection
  \ref{sec:hyperbolic} and subsection \ref{sec:proper}.
 
\item[(4)] \emph{Lagrange top}. The heavy top equations in body
  representation are known to be Hamiltonian on
  $\mathfrak{se}(3)^*$. These equations describe a classical
  Hamiltonian system with $2$ degrees of freedom on the coadjoint
  orbits of the Euclidean group $\op{SE}(3)$. The generic coadjoint
  orbit is a magnetic cotangent bundle of the 2-sphere $S^2$. This two
  degrees of freedom system has a conserved integral but it does not
  have, generically, additional integrals. However, in the Lagrange
  top, one can find one additional integral, namely, the momentum map
  associated with rotations about the axis of gravity, which makes the
  system integrable.  It is classically known that the Lagrange Heavy
  Top is integrable.
 
  The Lagrange momentum map is, however, a non\--proper map. But the
  singular Lagrangian fibration given by the system itself is a proper
  map. The Lagrange Top has hyperbolic singularities.
 
  Hence, this example fits in subsection \ref{sec:hyperbolic}
  (hyperbolic systems) and subsection \ref{sec:proper} (non\--proper
  systems). It may also overlap with subsection \ref{sec:degenerate}
  (degenerate systems), we have not checked this non\--degeneracy
  condition.

\item[(5)] \emph{Kowalevski top}.  This integrable system was discovered by
  Sophie Kowalevski and is published in her seminal paper \cite{Ko1889}.  This integrable system fits into a
  combination of subsection \ref{sec:antiperiodic} (non\--periodic
  systems) and subsection \ref{sec:degenerate} (degenerate systems).

  Although it is an integrable system, to our knowledge, \emph{it is
    not known} whether one of the components is periodic, i.e.,
  whether it comes from a Hamiltonian $S^1$\--action. It is generally
  believed that this is not the case but, as far as we know, there is
  no proof of this fact.

\item[(6)] \emph{Three wave interaction}. This example fits in subsection
  \ref{sec:higher} (higher\--dimensional systems), and subsection
  \ref{sec:degenerate} (degenerate systems). It may also overlap with
  subsection \ref{sec:hyperbolic} (hyperbolic systems). The
  phase\--space is $6$\--dimensional; we have checked that it has many
  degenerate singularities. The three wave interaction was brought up
  to our attention by D. Holm.

\end{itemize}

Examples (3), (4), (5), and (6) may not fit into \emph{any} of the
sections \ref{sec:hyperbolic}, \ref{sec:proper}, \ref{sec:higher},
\ref{sec:antiperiodic}, \ref{sec:degenerate}.  However, they do fit in
the program described in this paper, after the content of these
sections has been put together, as indicated in a subsection
\ref{sec:collaborative}.

\subsection{\textcolor{black}{Semitoric
    systems}} \label{sec:semitoric}

This section explains the technology, key tools, and recent techniques
developed to study the symplectic theory of integrable systems of
semitoric type. These methods have lead to a complete classification
in terms of five symplectic invariants \cite{PeVN2009, PeVN2011}.  The
overall strategy leading to these results should be applicable in a
much more general context. In this section we recall this strategy,
pointing out what technology particular to semitoric systems was used
and what methods are likely to extend to a more general context.

\subsubsection{Setting for integrable systems:
    symplectic manifolds} \label{sec:setting}
\begin{definition}
  A \emph{symplectic form} $\omega$ on a smooth manifold $M$ is a
  closed, non\--degenerate two\--form. The pair $(M, \, \omega)$ is
  called a \emph{symplectic manifold}.
\end{definition}

Let $G$ be a Lie group with Lie algebra $\mathfrak{g}$ whose dual is
denoted by $\mathfrak{g}^*$; $\exp:\mathfrak{g}\rightarrow G$ is the
exponential map.

\begin{figure}[h]
  \begin{center}
    \includegraphics[width=0.4\textwidth]{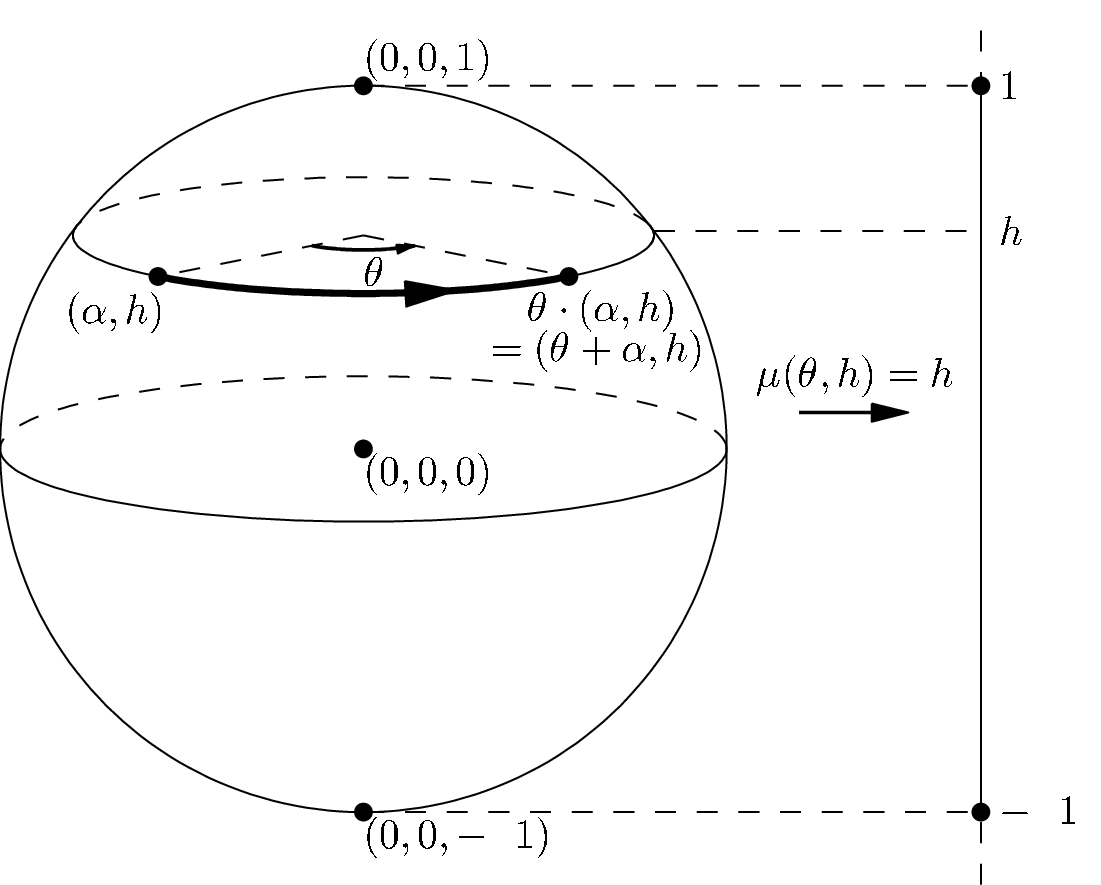}
    \caption{Momentum map $\mu(\theta,\, h)=h$ for the
      $2$\--dimensional sphere with the standard rotational
      $S^1$\--action. The momentum polytope in this case is
      $1$\--dimensional, the interval $[-1,\,1]$.}
    \label{fig:ff}
  \end{center}
\end{figure}

\begin{definition}
  A $G$\--action on $M$ is \emph{symplectic} if it preserves
  $\omega$. A symplectic $G$\--action is \emph{Hamiltonian} if there
  exists an equivariant map $\mu \colon M \to\mathfrak{g}^*$ such that
$$
\omega(X_M, \cdot)=  \op{d} \langle \mu, \, X\rangle,
$$
for all $X \in \mathfrak{g}$, where $X_M$ is the infinitesimal
generator vector field on $M$ induced by $X$, i.e., $$X_M (m) : =
\left.\frac{{\rm d}}{{\rm d}t}\right|_{t=0}\exp(tX) \cdot m$$ for
every $m \in M$; here $\langle \mu, \, X \rangle \colon M \to \R$ is a
smooth function obtained by using the duality pairing
$\left\langle\cdot,\, \cdot \right\rangle$ between $\mathfrak{g}^*$
and $\mathfrak{g}$.
\end{definition}

The map $\mu$ is called the \emph{momentum map} of the Hamiltonian
$G$\--action.  Since ${\rm i}_{X_M} \omega:=\omega(X_M,\, \cdot)$ is a
closed $1$\--form, every symplectic action is Hamiltonian if
$\op{H}^1(M,\, \R) = 0$. There are many examples of non\--Hamiltonian
$G$\--actions: take for instance the standard translational action of
the circle $S^1$ on the two\--torus $\mathbb{T}^2:=S^1 \times S^1$.

Our intuition on integrable systems has been guided by some remarkable
results proven in the 80s by Atiyah, Guillemin\--Sternberg and
Delzant, in the context of Hamiltonian torus actions.  Note that if
$\mathfrak{g}$ is $m$\--dimensional, by choosing a basis of
$\mathfrak{g}$ we may see the momentum map as a map into
$\mathbb{R}^m$ (see Figure \ref{fig:ff}).

\begin{theorem}[Atiyah \cite{atiyah} and Guillemin\--Sternberg \cite{gs}] \label{theo:ags} If an $m$\--dimensional torus
  $G$ with Lie algebra $\mathfrak{g}$ acts on a compact connected
  symplectic manifold $(M,\, \omega)$ in a Hamiltonian fashion, the
  image $\mu(M)$ of $M$ under the momentum map $ \mu \colon M \to
  \mathfrak{g}^* \simeq \R^m $ is a convex polytope.
\end{theorem}

\begin{figure}[htbp]
  \begin{center}
    \includegraphics[width=0.5\textwidth]{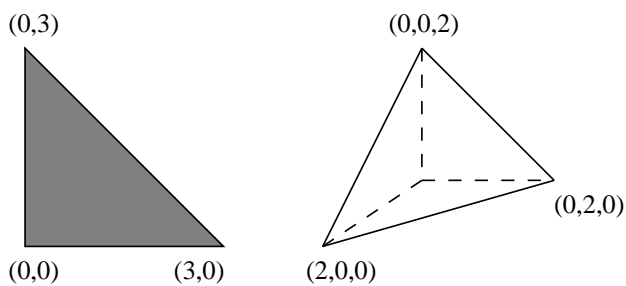}
    \caption{Delzant polytopes corresponding to the complex projective
      spaces $\C P^2$ and $\C P^3$ equipped with scalar multiples of
      the Fubini\--Study symplectic form.}
    \label{fig:polytopes}
  \end{center}
\end{figure}

Figure \ref{fig:polytopes} shows the images of the momentum map
described in Theorem \ref{theo:ags} for the standard toric actions on
complex projective spaces.

\begin{example}
  Consider the projective space $\mathbb{CP}^n$ equipped with a
  $\lambda$ multiple of the Fubini--Study form and the standard
  rotational action of $\mathbb{T}^n$ (for $\mathbb{CP}^1=S^2$, we
  already drew the momentum map in Figure \ref{fig:ff}). This action
  is Hamiltonian, and the momentum map is given by
$$
\mu(z)\,=(\frac{\lambda \,\, |z_1|^2}{\sum_{i=0}^n\,
  \,|z_i|^2},\ldots,\frac{\lambda \,\, |z_n|^2}{\sum_{i=0}^n\,
  \,|z_i|^2}).
$$
It follows that the image of $\mu$ equals the convex hull in
$\mathbb{R}^n$ of $0$ and the scaled canonical vectors $\lambda
e_1,\ldots,\lambda e_n$, see Figure \ref{fig:polytopes}.
\end{example}

Let ${\got g }_ \Z$ be the kernel of the exponential map
$\op{exp}:{\got g}\to G$. We denote the isomorphism ${\got g}/{\got
  g}_{\Z} \to G$ also by $\op{exp}$.
  
  \begin{definition}
    Let $G$ be an $n$\--dimensional torus, and let $\Delta \subset
    {\got g}^*$ be an $n$\--dimensional convex polytope.  Let
    $\mathcal{F}$ and $\mathcal{V}$ be the set of codimension one
    faces and vertices of $\Delta$, respectively. For $v\in
    \mathcal{V}$, write $$\mathcal{F}_v=\{ \ell \in \mathcal{F}\mid
    v\in f\}. $$ We say that $\Delta$ is a {\em Delzant polytope} in
    ${\got g}^*$ if:
    \begin{itemize}
    \item[i)] For each $\ell \in \mathcal{F}$ there exists $X_{\ell}
      \in {\got g}_{\Z}$ and $\lambda _{\ell} \in\R$ such that the
      hyperplane which contains $\ell$ is equal to the set of all
      $\xi\in {\got g}^*$ such that $\langle X_{\ell},\,\xi\rangle
      +\lambda _{\ell}=0$.
    \item[ii)] For every $v\in \mathcal{V}$, the vectors $X_{\ell}$
      with $\ell \in \mathcal{F}_v$ form a $\Z$\--basis of the
      integral lattice ${\got g}_{\Z}$ in ${\got g}$ (hence there are
      such $n$ vectors for each fixed $v$.)
    \end{itemize}
  \end{definition}

  Delzant showed the following stunning result.  In the statement of
  the theorem, the term \emph{isomorphism} between two symplectic
  manifolds $M$ and $M'$ respectively equipped with $G$\--actions for
  which the momentum maps are $\mu$ and $\mu'$, refers to a
  symplectomorphism\footnote{I.e., a diffeomorphism which pulls back
    the symplectic form on $M'$ to the symplectic form on $M$.}
  $\varphi \colon M \to M'$ such that $$\mu'\circ \varphi=\mu.$$ If it
  is the case, then $\varphi$ intertwines the torus actions.

  \begin{theorem}[Delzant \cite{delzant}] \label{delzantseminal} If the dimension $n$ of the torus $G$
    is half the dimension of $M$, then the polytope in the
    Atiyah\--Guillemin\--Sternberg Theorem is a Delzant polytope, this
    polytope determines the isomorphism type of $M$, and $M$ is a
    toric variety.  Even more, starting from any Delzant polytope in
    ${\got g}^*$, one can construct a symplectic manifold with a
    Hamiltonian $G$\--action for which its associated polytope is the
    one we started with.
  \end{theorem}

\begin{remark}
  Delzant's theorem says that the polytopes in Figure
  \ref{fig:polytopes} determine all the information about
  $\mathbb{CP}^n$, the symplectic form and the torus action on it.
\end{remark}

\begin{figure} [h]
  \centering
  \includegraphics[width=0.35\textwidth]{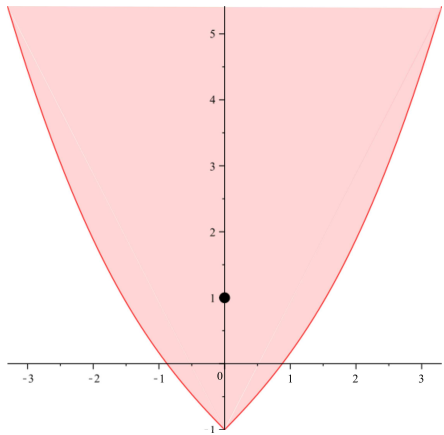}
  \caption{Image of the singular Lagrangian fibration given by the
    spherical pendulum.}
  \label{critical_set_spherical_pendulum.figure2}
\end{figure}

It is natural to wonder whether these striking results persist in the
case where the torus is replaced by a non\--compact group acting in a
Hamiltonian fashion on $M$. The seemingly simplest case happens when
the group is $\R^n$; the study of these $\R^n$\--actions is precisely
the goal of the theory of integrable systems.  The image of the
momentum map of an integrable system is usually not a convex polytope,
see Figure \ref{critical_set_spherical_pendulum.figure2}. In most
cases it is not even convex, and it would for instance have an annulus
shape.

The behavior of an integrable system -- which essentially always has
complicated singularities -- is much more flexible than that of a
torus action. Actions of $n$\--dimensional tori on $2n$\--dimensional
manifolds may be viewed as examples of integrable systems in the sense
that the components of the momentum map of the action form an
integrable system (we may also say the momentum map itself is an
integrable system).

\subsubsection{\textcolor{black}{Integrable systems}}

Let $(M,\omega)$ be a symplectic manifold of dimension $2n$.

\begin{definition}
  An \emph{integrable system} consists of $n$ functions $f_1, \ldots, f_n$ on
  $M$ whose differentials 
  $
  {\rm d}f_1, \ldots, {\rm d}f_n
  $
  are almost everywhere linearly independent
  and which Poisson-commute, i.e., $\{f_i,\,f_j\} = 0$ 
   for all integers $1 \leq i,\,j \leq n$.
   \end{definition}

We call $F:=(f_1,\ldots,f_n) \colon M \to \mathbb{R}^n$ the
\emph{momentum map} of the integrable system. Often we refer to the map 
$F$  simply as the \emph{integrable system}.

\begin{example}(Coupled spin\--oscillator) \label{ex:sp} A simple
  non\--toric non-compact example of integrable system is the so
  called ``coupled spin-oscillator'' model, which is of fundamental
  importance in physics, and therein known as the
  \emph{Jaynes\--Cummings model}.  In this case the symplectic
  manifold is $M=S^2\times\R^2$, where $S^2$ is viewed as the unit
  sphere in $\R^3$ with coordinates $(x,\,y,\,z)$, and the second
  factor $\R^2$ is equipped with coordinates $(u,\, v)$. The
  integrable system is given by the smooth maps
  $$
  J := (u^2+v^2)/2 + z$$ and $$H := \frac{1}{2} \, (ux+vy).$$
\end{example}

\begin{definition}
  A \emph{singularity}\footnote{We also call it  \emph{singular point} or \emph{critical point}.} of an integrable system 
  $$F:=(f_1,\ldots,f_n) \colon M \to \mathbb{R}^n$$
  is a point $m \in M$
  for which the tangent mapping ${\rm T}_mF: {\rm T}_mM \rightarrow
  \mathbb{R}^n$ has rank less than $n$, i.e.
${\rm d} f_1, \ldots, {\rm d}f_n$ are
linearly dependent one-forms at $m$.  The fiber that contains $m$ is
  a \emph{singular fiber}.
\end{definition}

\begin{figure} [h]
  \centering
  \includegraphics[width=\textwidth]{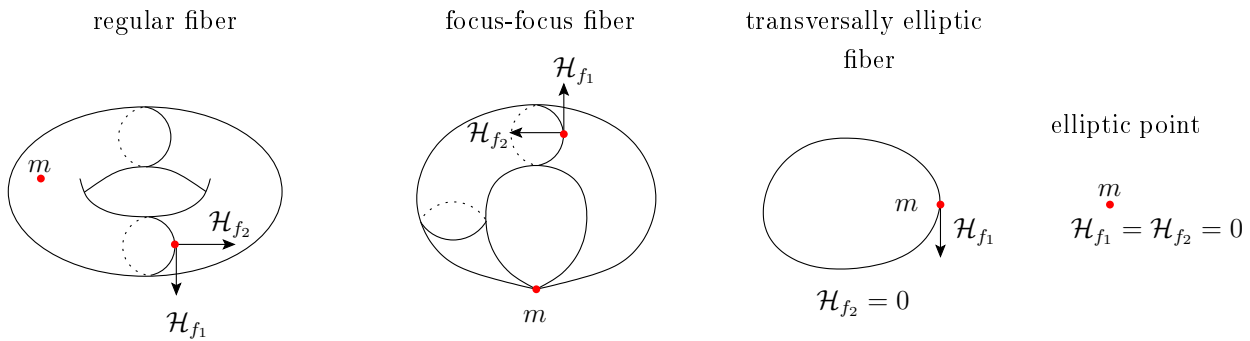}
  \caption{The figures show some possible singularities of an
    integrable system. On the left most figure, $m$ is a regular point
    (rank 2); on the second figure, $m$ is a focus-focus point (rank
    $0$); on the third one, $m$ is a transversally elliptic
    singularity (rank $1$); on the right most figure, $m$ is an
    elliptic-elliptic point (rank 0).}
  \label{singularities}
\end{figure}

The most interesting features of the integrable system are encoded in
the singular fibers of the momentum map, some of which are depicted in Figure \ref{singularities}.

\begin{definition}
  \label{def:non-deg}
  Let $F=(f_1,\dots,f_n) \colon M \to \mathbb{R}^n$ be an integrable
  system and let $m\in M$ be a singularity.
  \begin{itemize}
\item[(a)] If the corank of the tangent map ${\rm T}_mF: {\rm T}_m M
  \rightarrow \mathbb{R}^n$ is $n$, the singularity $m$ is said to
  be \emph{non-degenerate} if the Hessians $\op{d}^2_mf_j$ span a
  Cartan subalgebra of the Lie algebra of quadratic forms on the
  symplectic vector space $({\rm T}_m M, \omega_m)$.
\item[(b)] If the corank of 
  ${\rm T}_mF: {\rm T}_m M \rightarrow \mathbb{R}^n$
is $r=0, \ldots , n-1$, we can assume without loss of generality that
$
{\rm d}_mf_1,\dots,{\rm d}_mf_{n-r}
$
  are linearly independent.  Let $\Sigma$ be a $2r$-dimensional local
  symplectic submanifold containing $m$ and transversal to the flows
  of the Hamiltonian vector fields $\mathcal{H}_{f_1},\ldots, \mathcal{H}_{f_{n-r}}$ at $m$ (defined
  by $\omega(\mathcal{H}_{f_i},\,\cdot)={\rm d}f_i$). Then
  the rank of $$(f_{n-r+1}|_\Sigma, \dots, f_n|_ \Sigma): \Sigma
  \rightarrow \mathbb{R}^{n-r}$$ at $m$ is zero. We say that $m$ is
  \textit{non-degenerate}, or \textit{transversally non-degenerate} at
  $m$ if this new integrable system on $\Sigma$ has a rank zero
  non-degenerate singularity (in the sense defined above).
\end{itemize}
\end{definition}

It follows from the work Williamson in the 1930s that a Cartan
subalgebra as in Definition \ref{def:non-deg} has a basis with three
type of blocks:
\begin{itemize}
\item[(1)] one-dimensional ones:
  \begin{itemize}
  \item \emph{elliptic block} $q=x^2+\xi^2$,
  \item \emph{hyperbolic block} $q=x\xi$;
  \end{itemize}
\item[(2)] two-dimensional one: \emph{focus-focus block}:
  $q_1=x\eta-y\xi$, $q_2=x\xi+y\eta$.
\end{itemize}

\begin{example}
  The coupled spin\--oscillator system in Example \ref{ex:sp} has
  non\--degenerate singularities only. It has exactly one
  focus\--focus singularity at the ``North Pole''
  \[
  ((0,\,0,\,1),\,(0,\,0)) \in S^2 \times \R^2
  \]
  and one elliptic\--elliptic singularity at the ``South Pole''
  $((0,\,0,\,-1),\,(0,\,0))$. The remaining singularities are
  transversally elliptic.
\end{example}

The following notion is due to N.T. Zung.

\begin{definition} \label{def:zung} Let $m$ be a singularity of an
  integrable system $F \colon M \to \mathbb{R}^n.$  If $k_{\op{e}},\,
  k_{\op{h}},\, k_{\op{f}}$ respectively denote the number of
  elliptic, hyperbolic, and focus\--focus components of $m$, we call
  $(k_{\op{e}},\, k_{\op{h}},\, k_{\op{f}})$ the \emph{Williamson type
    of $m$}.
\end{definition}

The following theorem is one of the key tools in the modern theory of
integrable systems. It was priorly proven in the
analytic case by Vey.

\begin{theorem}[Eliasson \cite{eliasson}] \label{singularities_theorem} The
  non\--degenerate singularities of an integrable system
  $F \colon M \to \mathbb{R}^n$ are linearizable, i.e., if $m \in M$
  is a non\--degenerate singularity of the completely integrable
  system $F=(f_1,\,\ldots,f_n): M \rightarrow \mathbb{R}^n$ then there
  exist local symplectic coordinates $(x_1,\, \ldots,x_n,\, \xi_1,\,
  \ldots,\, \xi_n)$ about $m$, in which $m$ is represented as
  $(0,\,\ldots,\, 0)$, such that $\{f_i,\,q_j\}=0$, for all indices
  $i,\,j$, where we have the following possibilities for the
  components $q_1,\,\ldots,\,q_n$, each of which is defined on a small
  neighborhood of $(0,\,\ldots,\,0)$ in $\mathbb{R}^n$:
  \begin{itemize}
  \item[{\rm (i)}] Elliptic component: $q_j = (x_j^2 + \xi_j^2)/2$,
    where $j$ may take any value $1 \le j \le n$.
  \item[{\rm (ii)}] Hyperbolic component: $q_j = x_j \xi_j$, where $j$
    may take any value $1 \le j \le n$.
  \item[{\rm (iii)}] Focus\--focus component: $q_{j-1}=x_{j-1}\,
    \xi_{j} - x_{j}\, \xi_{j-1}$ and $q_{j} =x_{j-1}\, \xi_{j-1}
    +x_{j}\, \xi_{j}$ where $j$ may take any value $2 \le j \le n-1$
    (note that this component appears as ``pairs'').
  \item[{\rm (iv)}] Non\--singular component: $q_{j} = \xi_{j}$, where
    $j$ may take any value $1 \le j \le n$.
  \end{itemize}
  Moreover if $m$ does not have any hyperbolic component, then the
  system of commuting equations $\{f_i,\,q_j\}=0$, for all indices
  $i,\,j$, may be replaced by the single equation
 $$
 (F-F(m))\circ \varphi = g \circ (q_1,\, \ldots,\,q_n),
 $$
 where $\varphi=(x_1,\, \ldots,x_n,\, \xi_1,\, \ldots,\, \xi_n)^{-1}$
 and $g$ is a diffeomorphism from a small neighborhood of the origin
 in $\mathbb{R}^n$ into another such neighborhood, such that $g(0,\,
 \ldots,\,0)=(0,\,\ldots,\,0)$.
\end{theorem}

If the dimension of $M$ is $4$ and $F$ has no hyperbolic
singularities, we have the following possibilities for the map
$(q_1,\,q_2)$, depending on the rank of the singularity:
\begin{itemize}
\item[{\rm (1)}] if $m$ is a singularity of $F$ of rank zero, then
  $q_j$ is one of
  \begin{itemize}
  \item[{\rm (i)}] $q_1 = (x_1^2 + \xi_1^2)/2$ and $q_2 = (x_2^2 +
    \xi_2^2)/2$.
  \item[{\rm (ii)}] $q_1=x_1\xi_2 - x_2\xi_1$ and $q_2 =x_1\xi_1
    +x_2\xi_2$; \,\, on the other hand,
  \end{itemize}
\item[(2)] if $m$ is a singularity of $F$ of rank one, then
  \begin{itemize}
  \item[{\rm (iii)}] $q_1 = (x_1^2 + \xi_1^2)/2$ and $q_2 = \xi_2$.
  \end{itemize}
\end{itemize}

\begin{definition}
  Suppose that $(M,\omega)$ is a $4$\--dimensional symplectic manifold 
  and  that $F \colon M \to
  \mathbb{R}^2$ is an integrable system. A non\--degenerate singularity 
  is respectively called \emph{elliptic\--elliptic,
    focus\--focus}, or \emph{transversally\--elliptic} if both
  components $q_1,\, q_2$ above are of elliptic type, $q_1,\,q_2$
  together correspond to a focus\--focus component, or one component
  is of elliptic type and the other component is $\xi_1$ or $\xi_2$,
  respectively.   Similar definitions hold for \emph{transversally-hyperbolic,
    hyperbolic-elliptic} and \emph{hyperbolic\--hyperbolic}
  non-degenerate singularity.
\end{definition}

\subsubsection{\textcolor{black}{Semitoric systems: definition and
    basics}}

The classification theory we outline next was introduced by the
authors in \cite{PeVN2009,PeVN2010}.  The theory of semitoric 
integrable systems is at the early development
stages, which have already yielded results beyond our initial
expectations.  Our estimate is that it may take many years of work to
push the theory (both at the symplectic and spectral level) to cover
integrable systems which are \emph{not} necessarily
semitoric\footnote{Here we are including integrable systems on higher
  dimensional manifolds. Though semitoric systems may naturally be
  defined in all dimensions, they were originally defined for
  $4$\--dimensional manifolds.}.

We start with the basic notions. \emph{For the remaining of Section
  \ref{sec:semitoric} we work only with connected four-dimensional
  manifolds.}

\begin{definition} \label{def:semitoric} Let $(M,\, \omega)$ be a
  connected symplectic manifold.  An integrable system $F:=(J,\,H)
  \colon M \to \mathbb{R}^2$ is \emph{semitoric} if $J$ is a proper
  momentum map for an effective $S^1$ action on $M$ and $F$ has only
  non-degenerate singularities without hyperbolic singularity.
\end{definition}

Let us spell out Definition \ref{def:semitoric} more concretely. A
\emph{semitoric system} consists of a connected symplectic
four\--dimensional manifold $(M,\, \omega)$ and two smooth functions
$J \colon M \to \R$ and $H \colon M \to \R$ such that:
\begin{itemize}
\item[(a)] $J$ is constant along the flow of the Hamiltonian vector
  field $\mathcal{H}_H$ generated by $H$ or, equivalently,
  $\{J,\,H\}=0$;
\item[(b)] for almost all points $p \in M$, the vectors
  $\mathcal{H}_J(p)$ and $\mathcal{H}_H(p)$ are linearly independent;
\item[(c)] $J$ generates a $2\pi$\--periodic flow, i.e., $J$ is the
  momentum map of an $S^1$\--action on $M$;
\item[(d)] $J$ is a proper map;
\item[(e)] $F$ has only non\--degenerate singularities without
  hyperbolic components.
\end{itemize}

\begin{remark}
  \normalfont A semitoric system has only elliptic and focus-focus
  singularities.
\end{remark}

\begin{definition}\label{iso_def}
  Two semitoric systems 
  $$(M_1,\omega_1,(J_1,H_1)) \,\, \textup{and}\,\, (M_2,\omega_2,(J_2,H_2)) 
  $$
  are
  \emph{isomorphic} if there exists a symplectomorphism $\phi \colon
  M_1\to M_2$, and a smooth map $f \colon F_1(M_1)\to \mathbb{R}$ with
  $\partial_2f \neq 0$, such that $J_1 = \phi^*J_2$ and $\phi^*H_2 =
  f(J_1,H_1)$.
\end{definition}

``Smooth" for the map $f$ means that it is the restriction of a smooth
map $U \rightarrow \mathbb{R}$, where $U$ is an open subset of
$\mathbb{R}^2$ containing $F_1(M_1)$.

\begin{example}
  Perhaps the simplest non\--toric, semitoric integrable system on a
  non\--compact manifold is the coupled spin\--oscillator model
  described Example \ref{ex:sp}.
\end{example}

 A polygon $\Delta
\subset \mathbb{R}^2$ is said to be \textit{rational}, if the normal
vectors to its edges span a sublattice of $\mathbb{Z}^2$.

\begin{theorem}[V\~ u Ng\d oc \cite{VN2007}] \label{advances} To a
  semitoric integrable system $F =(J, H) \colon M \to \mathbb{R}^2$ one can
  associate a convex rational polygon $\Delta \subset \mathbb{R}^2$,
  unique up to translations, and modulo the action of the matrix group
  with elements
$$
\left(
  \begin{array}{cc}
    1 & 0\\ k & 1
  \end{array}
\right).
  $$
  In addition, $\Delta= \mu(F(M))$, where $\mu$ is a homeomorphism and
  $\mu$ is integral-affine\footnote{For the exact technical definition
    of the notion of integral-affine, see the review \cite{PeVN2012}.}
  on a dense open subset.
\end{theorem}

\begin{figure}[htb]
  \begin{center}
    \includegraphics[width=0.6\textwidth]{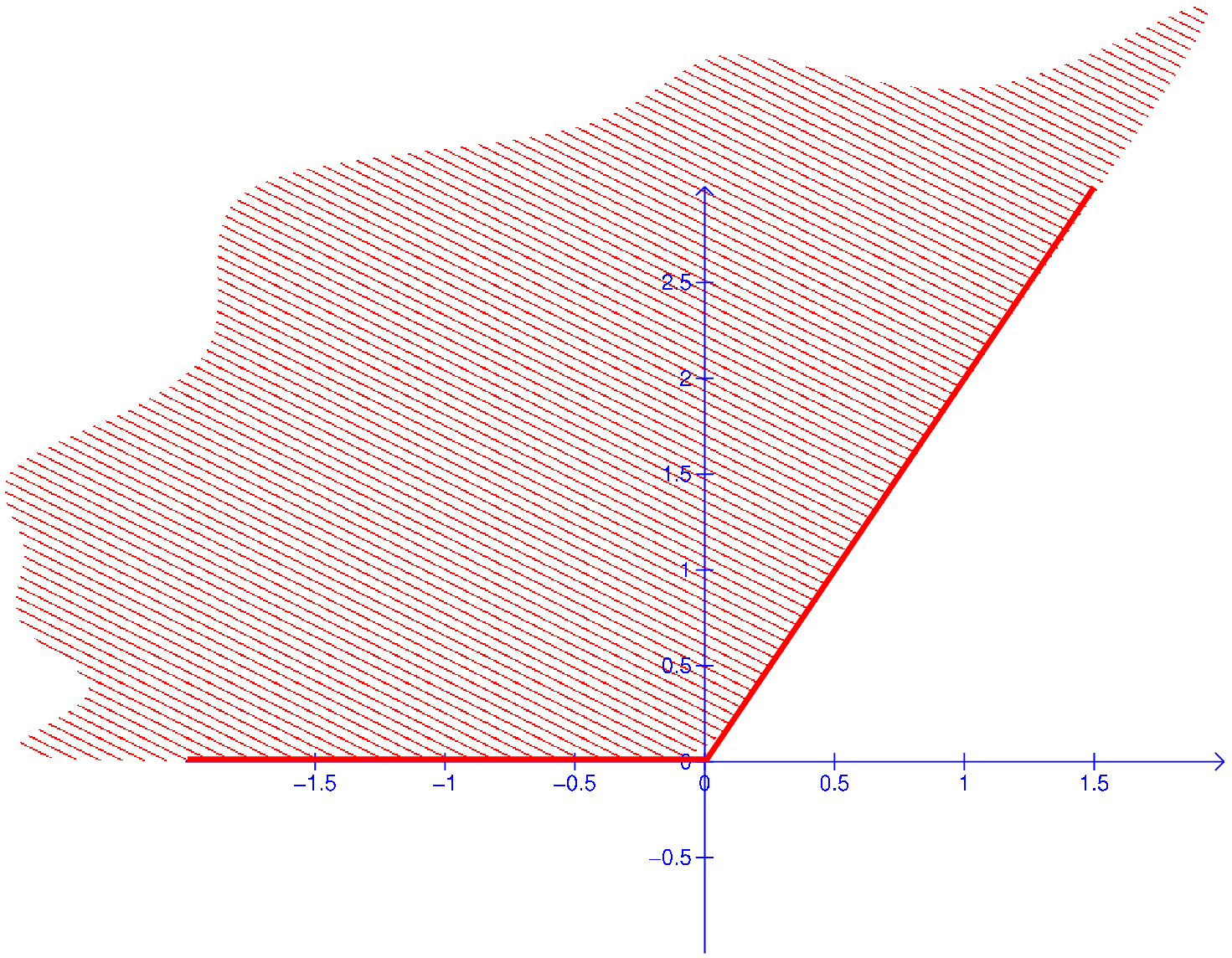}
    \caption{Possible transformed image $\mu(F(M))$ of a semitoric
      system $F \colon M \to \mathbb{R}^2$, see Theorem
      \ref{advances}.}
  \end{center}
\end{figure}

It is natural to wonder whether the polygon constructed in Theorem
\ref{advances} characterizes the system, as in the Delzant's
classification result (Theorem \ref{delzantseminal}).  We elaborate on
this more in subsection \ref{sec:higher}.

\subsubsection{\textcolor{black}{Uniqueness theorem for semitoric
    systems}}

The polygon in Theorem \ref{advances} does not characterize the
system.  However, a larger collection of invariants, which includes
the polygon, does determine the system. This is the content of the
following result.

\begin{theorem}[Pelayo-V\~ u Ng\d oc \cite{PeVN2009}] \label{inventiones} A semitoric system is
  characterized, up to isomorphisms, by the following invariants:
  \begin{itemize}
  \item[{\rm (1)}] {\bf \emph{number of singularities}}: the number of
    focus\--focus singularities $m_f$;
  \item[{\rm (2)}] {\bf \emph{singular foliation type}}: a formal
    Taylor series $S(X,Y)$ at each focus\--focus singularity;
  \item[{\rm (3)}] {\bf \emph{polygon invariant}}: a class of polygons
    equipped with $m_f$ oriented vertical lines (see Figure
    \ref{fig:weightedpolygon});
  \item[{\rm (4)}] {\bf \emph{volume}}: $m_f$ points
    $c_1,\ldots,c_{m_f}$ contained in the interior of the image of the
    system or, equivalently, a finite set of positive numbers giving
    the positions of these interior points;
  \item[{\rm (5)}] {\bf \emph{twisting\--index}}: a collections of
    $m_f$ integer (one integer between two consecutive nodes ordered
    according to their $J$\--component).
  \end{itemize}
\end{theorem}

In other words, Theorem \ref{inventiones} says that if $(M,\,
\omega_1,\,(J_1,\,H_1))$ and $(M,\, \omega_2,\, (J_2,\,H_2))$ are
semitoric systems, then
\begin{center}
  \emph{$(M,\, \omega_1,\,(J_1,\,H_1))$ and $(M,\, \omega_2,\,
    (J_2,\,H_2))$ are isomorphic $\iff$ their invariants (1)--(5)
    agree.}
\end{center}
The word \emph{isomorphism} is used in the sense of Definition
\ref{iso_def}, i.e., there exists a symplectomorphism
$$
\varphi \colon M_1 \to M_2,\,\,\,\,\,\,\, \textup{such that}\,\,\,\,\,
\varphi^*(J_2,\,H_2)=(J_1,\,f(J_1,\,H_1)).
$$
for some smooth function $f \colon F_1(M_1)\to \mathbb{R}$ with
  $\partial_2f \neq 0$.

Let's comment on the symplectic invariants above; their precise
construction is found in \cite[Section 2]{PeVN2009}.
\begin{itemize}
\item[(1)] \emph{The number of singularities} is an integer $m_f$
  counting the number of isolated singularities of the integrable
  system (which correspond precisely to the images of singularities of
  focus-focus type).  We write $m_f$ to emphasize that the
  singularities that $m_f$ counts are focus\--focus singularities.
\item [(2)] \emph{The singular foliation type} is a collection of
  $m_f$ infinite Taylor series on two variables which classifies
  symplectically a saturated neighborhood of the singular fiber.
\item[(3)] \emph{The polygon invariant} is the equivalence class of a
  weighted rational convex\footnote{generalizing the Delzant polygon
    and which may be viewed as a bifurcation diagram} polygon
$$
\left(\Delta,\,(\ell_j)_{j=1}^{m_f},\,(\epsilon_j)_{j=1}^{m_f}\right),
$$
which is constructed from the image of the system by performing a very
precise ``cutting" which appears in the proof of Theorem
\ref{advances} (see Figure \ref{fig:weightedpolygon}).  Here $\Delta$
is a convex polygonal domain in $\R^2$, the $\ell_j$ are vertical
lines intersecting $\Delta$ and the $\epsilon_j$ are $\pm 1$ signs
giving each line $\ell_j$ an orientation.
\begin{figure}[htb]
  \begin{center}
    \includegraphics[width=0.5\textwidth]{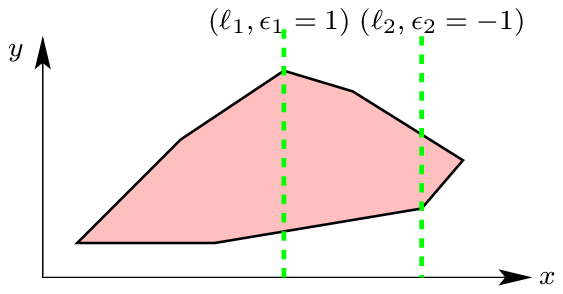}
    \caption{Weighted polygon.}
    \label{fig:weightedpolygon}
  \end{center}
\end{figure}
\item[(4)] \emph{The volume} invariant consists of $m_f$ real numbers
  giving the volumes of certain submanifolds meeting at the
  singularities; this invariant can be characterized by the position
  of a number of interior points $c_j$ in $\Delta$, which correspond
  to the focus\--focus values of the system.
\item[(5)] \emph{The twisting index} consists of $m_f$ integers
  measuring how twisted the Lagrangian fibration is around the
  singularities. This is a subtle invariant, which depends on the
  representative chosen in (3). \end{itemize}

\begin{remark}
  In Theorem \ref{inventiones}, the ``polygon invariant'' encodes the
  regular points and the elliptic singularities. The other invariants
  are linked to focus-focus singularities.
\end{remark}

\textup{\,}
\\
{\emph{Proof Strategy for Theorem \ref{inventiones}.}
The proof strategy for the uniqueness part is simple. We   start with two
integrable systems $F_1=(J_1,\,H_1) \colon M\to \mathbb{R}^2$ and
$F_2=(J_2,\,H_2) \colon M \to \mathbb{R}^2$.
\\
\\
\emph{Step 1} (Construction of invariants).  Construct local
symplectic invariants which are analogues of the invariants (1)--(5)
above.
\\
\\
\emph{Step 2} (Construction of local symplectomorphisms). Argue that
one can reduce to a case where the images $F_1(M_1)$ and $F_2(M_2)$
are equal. Prove that this common image can be covered by open sets
$\Omega_\alpha$, above each of which $F_1$ and $F_2$ are
symplectically interwined, i.e., above each $\Omega_{\alpha}$ one
finds a part of the system $F_1$ and a part of the system $F_2$ that
are symplectically the same because they have the same symplectic
invariants. At this stage one has a collection of local
symplectomorphisms which cover all the pieces of the manifolds $M_1$
and $M_2$.
\\
\\
\emph{Step 3} (Symplectic gluing). Use symplectic gluing of these
local pieces to prove a uniqueness and an existence type theorem. The
last step is to glue together these local symplectomorphisms, thus
constructing a global isomorphism $\phi:M_1\to M_2$.

\vspace{1mm}

\emph{We expect that the same proof strategy applies to each of the
  upcoming subsections, so we do not repeat it therein. Instead, each
  of the following subsections is focused on explaining what the
  difficulties are, mainly in Step 1 (but we also comment on potential
  difficulties in Step 2 and Step 3).}

\begin{remark}
  The ``analytic-combinatorial" data in Theorem \ref{inventiones}
  completely describes the moduli space of semitoric systems, as shown
  in the following reconstruction theorem. These invariants are
  depicted in Figure \ref{fig:spin} for the case of the coupled
  spin\--oscillator.
\end{remark}

  \begin{figure}[h]
    \centering
    \includegraphics[width=0.9\linewidth]{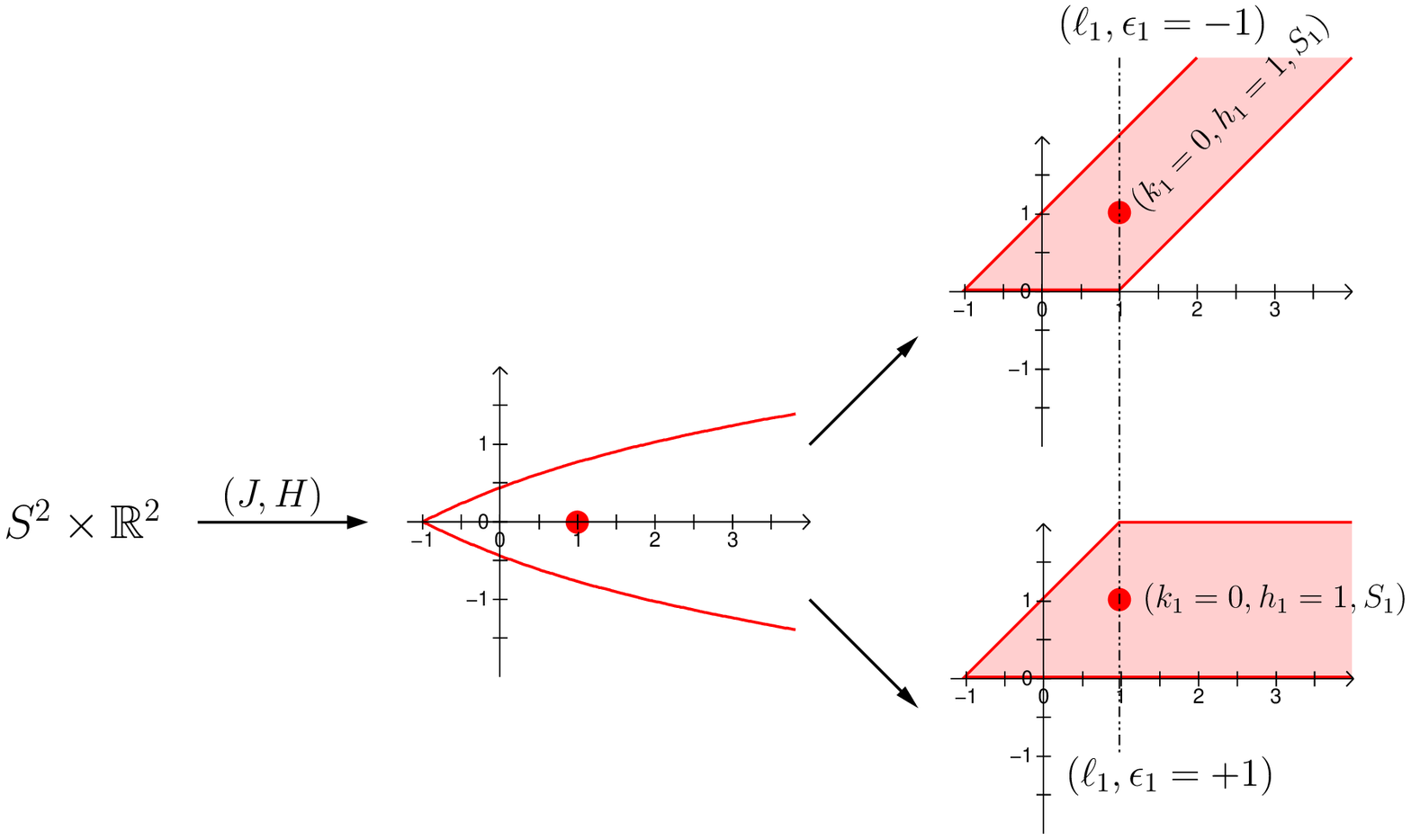}
    \caption{The coupled spin-oscillator example is a semitoric system
      with one focus\--focus point whose image is $(1,\,0)$. The
      invariants are depicted on the right hand-side.  In this case
      $m_f=1$ and the class of generalized polygons for this system
      consists of two polygons. $S_1$ denotes the Taylor series
      invariant (the linear terms of which were computed in
      \cite{PeVN2010}), $k_1=0$ is the twisting index invariant (which
      is trivial because there is only one focus\--focus value in this
      case), and $h=1$ is the volume invariant, which determines the
      position of the focus\--focus value node.}
    \label{fig:spin}
  \end{figure}

  \subsubsection{\textcolor{black}{Existence theorem for semitoric
      systems}}

  The following result appeared in \cite[Theorem 4.7]{PeVN2011}.

  \begin{theorem}[Pelayo-V\~ u Ng\d oc \cite{PeVN2011}]
    \label{acta}
    Given the following ingredients:
    \begin{itemize}
    \item[{\rm (1)}] {\bf \emph{number of singularities}}: an integer
      number $0 \le m_f<\infty$;
    \item[{\rm (2)}] {\bf \emph{singular foliation type}}: an
      $m_{f}$\--tuple of Taylor series \[((S_i)^{\infty})_{i=1}^{m_{f}}
      \in (\R[[X,\,Y]]_ 0)^{m_f};\]
    \item[{\rm (3})] {\bf \emph{polygon}}: a Delzant semitoric polygon
      $[\Delta_{\scriptop{w}}]$ of complexity $m_f$ (see {\rm
        \cite[\emph{Definition 4.3} ]{PeVN2010} for the notion of
        Delzant semitoric polygon}).  We denote the representative
      $\Delta_{\scriptop{w}}$ of $[\Delta_{\scriptop{w}}]$ by
      \[\Big(\Delta,\, (\ell_{\lambda_j})_{j=1}^{m_f},\,
      (\epsilon_j)_{j=1}^{m_f}\Big);\]
      
    \item[{\rm (4)}] {\bf \emph{volume}}: an $m_f$\--tuple of numbers
      $(h_j)_{j=1}^{m_{f}}$ such that \[0 < h_j <
      \textup{length}(\Delta \cap \ell_i)\] for each $j \in
      \{1,\ldots,m_f\}$;
    \item[{\rm (5)}] {\bf \emph{twisting\--index}}: a collection of
      $m_f-1$ integers (up to some equivalence relation, see {\rm
        \cite[\emph{Theorem 4.7}]{PeVN2011}}).
    \end{itemize}
    Then there exists a symplectic $4$-manifold $(M,\, \omega)$ and a
    semitoric integrable system $F = (J,\,H)$ on $M$ whose symplectic
    invariants of Theorem \ref{inventiones} given in that order
    coincide with (1)--(5) above.
  \end{theorem}

\textup{\,}
\\
{\emph{Proof Strategy for Theorem \ref{acta}.}}
The proof strategy for the existence theorem is simple.
\\
\\
\emph{Step 1} (Start with abstract collection of items
``corresponding" to invariants).
Let $$\Big(\Delta,\,(\ell_{\lambda_j})_{j=1}^s,\,
(\epsilon_j)_{j=1}^s\Big)$$ be a representative of the polygon
invariant with all $\epsilon_{j}$'s equal to $+1$. The strategy is to
use a symplectic gluing procedure introduced in \cite{PeVN2011} in
order to obtain a semitoric system by constructing a suitable singular
torus fibration above $\Delta\subset\RM^2$.

For $j=1,\dots, m_f$, let $c_j\in\RM^2$ be the point with coordinates
\begin{equation}
  c_j=(\lambda_j,\,h_j+\min(\pi_2(\Delta\cap\ell_{\lambda_j}))),
  \label{equ:cuts} \nonumber
\end{equation}
where $\pi_2: \mathbb{R}^2\rightarrow \mathbb{R}$ is the projection on
the second factor.  Because of the assumption on $h_j$, all points
$c_j$ lie in the interior of the polygon $\Delta$. We call these
points \emph{nodes}.  We denote by $\ell_j^+$ the vertical half\--line
through $c_j$ pointing upwards. We call these half-lines \emph{cuts}.
Now construct a ``convenient" (this is technical, we leave it to the
papers) covering of the polygon $\Delta$.
\\
\\
\emph{Step 2} (Local construction piece by piece). Construct a
``semitoric system'' over the part of the polygon away from the sets
in the covering that contain the cuts $\ell_j^+$; then we attach to
this ``semitoric system'' the focus\--focus fibrations, i.e., the
models for the systems in a small neighborhood of the nodes.  Continue
to glue the local models in a small neighborhood of the cuts. The
``semitoric system'' is given by a proper toric map only in the
preimage of the polygon away from the cuts.
\\
\\
\emph{Step 3} (Recovering smoothness).  Modify the system to recover
the smoothness of the system (this is very delicate) and observe that
the invariants of the system are precisely the items (1)--(5) we
started with.

\subsection{\textcolor{black}{Hyperbolic
    systems}} \label{sec:hyperbolic}

This branch of the program has been explored by Zhilinski\'{\i} from
the point of view of physics.  However it is not rigorously formalized
as a mathematical theory. Zhilinski\'{\i}'s explorations are
intriguing, and should give a lot of insight when working towards a
mathematical theory. There are also a few related results by V\~{u}
Ng\d{o}c, Dullin Zung, and Bolsinov, among others. The
one\--dimensional case (i.e. when the phase space is a symplectic
$2$\--dimensional manifold) was solved by Toulet, Molino, and Dufour.

Next we elaborate on the difficulties to construct the symplectic
invariants for this type of integrable systems. Once this is done, the
goal is to achieve a classification (existence and uniqueness) result
in terms of these invariants.

\begin{figure}[htbp]
  \begin{center}
    \includegraphics[width=0.8\textwidth]{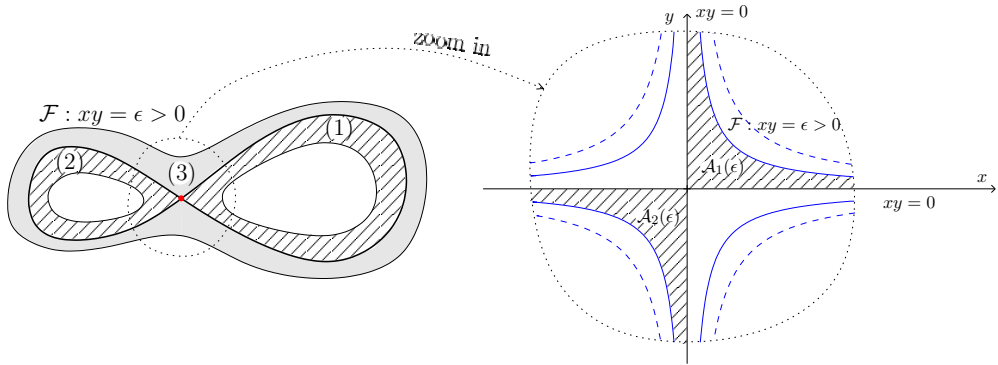}
    \caption{Zoom in around a \emph{hyperbolic component} of a
      singularity.}
    \label{hyperbolic-zoom}
  \end{center}
\end{figure}

\emph{We propose a general strategy to prove such classification
  immediately after the statements of Theorem \ref{inventiones} and
  Theorem \ref{acta} in the previous sections, so we do not repeat it
  here.} Instead, we explain what the expected difficulties are, as
far as we can see, to construct the invariants in the case of the
systems treated in this section. \emph{This comment also applies to
  the upcoming subsections \ref{sec:proper}, \ref{sec:higher},
  \ref{sec:antiperiodic}, \ref{sec:degenerate} and
  \ref{sec:topology}.}

In what follows, we refer to the previous subsection for the basic
terminology and results (and to \cite[Sections 3--8]{PeVN2012} and the
references therein for more detailed explanations).

\subsubsection{Assumptions} We say that a non-degenerate singularity
of an integrable system is of \emph{hyperbolic type} if the Willamson
type of $m$ (see Definition \ref{def:zung}) has $k_{\op{h}} \neq 0$.
In other words, a non-degenerate singularity has hyperbolic type if it
has some hyperbolic component.  Suppose that $F \colon M \to
\mathbb{R}^2$ is an integrable system which fails to be semitoric only
because it has some singularities containing components of hyperbolic
type. Concretely, this means that $(M,\, \omega)$ is a connected
symplectic $4$\--manifold equipped with two smooth functions $J \colon
M \to \R$ and $H \colon M \to \R$ such that:
\begin{itemize}
\item[(a)] $J$ is constant along the flow of the Hamiltonian vector
  field $\mathcal{H}_H$ generated by $H$ or, equivalently,
  $\{J,\,H\}=0$;
\item[(b)] for almost all points $p \in M$, the vectors
  $\mathcal{H}_J(p)$ and $\mathcal{H}_H(p)$ are linearly independent;
\item[(c)] $J$ generates a $2\pi$\--periodic flow, i.e., $J$ is the
  momentum map of an $S^1$\--action on $M$;
\item[(d)] $J$ is a proper map;
\item[(e)] $F$ has only non\--degenerate singularities (possibly with
  hyperbolic singularities).
\end{itemize}

The following is the classification problem of the section.

\begin{classproblem} \label{problemv1} \normalfont The problem has
  three parts.
  \begin{itemize}
  \item[{\rm (I)}] Give explicit constructions of ``analogues" of the
    symplectic invariants in Theorem \ref{inventiones} (which refer to
    semitoric systems only), for integrable systems satisfying
    assumptions (a)\--(e).
  \item[{\rm (II)}] Define an abstract list of all such possible
    invariants which occur in (I).
  \item[{\rm (III)}] Extend Theorems \ref{inventiones}, \ref{acta} to
    integrable systems satisfying (a)\--(e), using (I) and (II).
  \end{itemize}
\end{classproblem}

At this time we are not aware of examples of systems satisfying all
assumptions (a)\--(e) but we do not see any a priori reason why they
would not exist. In fact, it is through the study of invariants that
we propose in this section that one should be able to construct large
classes of them, if in fact they exist (this was the case in our
article \cite{PeVN2011}), or disprove their existence. Nevertheless
there are physically relevant examples like the two\--body problem or
the Lagrange Top that satisfy several of these assumptions, see
Section \ref{sec:first}.

\begin{figure}[h]
  \centering
  \includegraphics[width=0.5\linewidth]{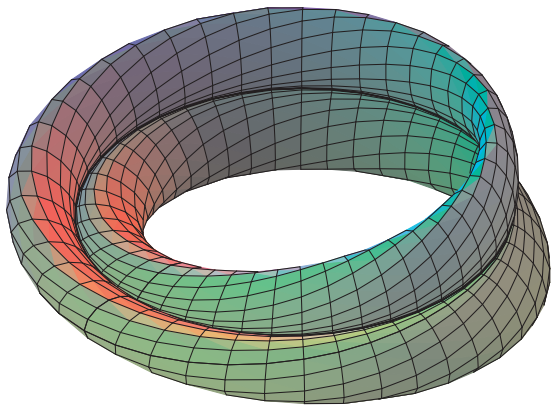}
  \caption{Hyperbolic critical fiber with $\Z_2$ symmetry.}
  \label{fig:hypspectra}
\end{figure}

\subsubsection{Expected difficulties} The expected difficulties are:
\begin{itemize}
\item[(i)] {\bf \emph{Polygon invariant}}.  There is no obvious
  candidate for the polygon invariant.

  The original construction of the polygon invariant from the image of
  $F(M)$ is by cutting along the vertical lines that pass through the
  focus\--focus values of $F$. At an ideological level, the
  construction works because focus\--focus singularities are isolated
  and the behavior of the integrable system around a focus\--focus
  singularity (while very complicated and having infinitely many
  symplectic invariants) can be controlled in a small semiglobal
  neighborhood $\mathcal{U}$ of the focus\--focus singularity. In
  particular, we know that focus\--focus singular fibers of semitoric
  systems are connected. However, singular fibers over hyperbolic
  points may not be connected, and the number of connected components
  changes when passing through it. This fiber connectivity is
  essential in the construction of the polygon invariant for semitoric
  systems\footnote{it is also intimately related to the convexity of
    the image of the momentum map in the
    Atiyah\--Guillemin\--Sternberg theorem, see Theorem
    \ref{theo:ags}.}.

  The image of $\mathcal{U}$ corresponds then to a vertical band in
  the image $F(M)$, which contains the vertical line passing through
  the focus\--focus singularity. This is, however, not the case with
  hyperbolic singularities, where the singularity may, for instance,
  come as a curve of singularities (e.g., suppose that the singularity
  has a hyperbolic component and a regular component) and it is not
  clear what the image of this curve is inside of $M$, and its
  relation to the integral affine structure induced by the the system.

  Understanding what the image of the hyperbolic singularity and its
  surrounding singularities are is a must if one wants to construct
  any kind of polygon invariant. Once this is clarified, it is not
  clear whether an invariant as simple as the polygon invariant may be
  constructed. It is more likely that a new more complicated invariant
  is going to replace the polygon invariant, probably a
  \emph{foliation type invariant}: a collection of $2$\--dimensional
  leaves with some structure. Foliations, and foliation type
  invariants, are important in many contexts, see for instance
  Bolsinov\--Fomenko \cite{bolsinov-fomenko-book} from a singularity
  theory angle and Bramham\--Hofer \cite[Section 4]{[14]} from the point
  of view of holomorphic curves.

\item[(ii)] {\bf \emph{Other invariants}}.  The definition of the
  other symplectic invariants should be less difficult than (i).
  However, constructing the twisting\--index invariant in the presence
  of hyperbolic singularities is unclear to us. The original
  construction of this invariant depended on the polygon invariant,
  which is as yet not known. Also, the twisting index\--invariant
  relied on an ordering of the focus\--focus singularities. If there
  are hyperbolic singularities, it is possible that one needs to
  construct a new twisting-index invariant.

  The original construction of the twisting\--index invariant encodes
  the ``difference'' between consecutive (according to the ordered
  $J$\--values) normalizations of the integrable system $F$, which is
  given by a $2$ by $2$ integer matrix with all entries constant, with
  the exception of one entry which, in general, is a non\--zero
  integer. This non\--zero integer is precisely the
  twisting\--index. This matrix is very much linked to the nature of
  the focus\--focus singularities.

  Because we are in dimension four, and hence focus\--focus
  singularities are isolated, it is quite possible that we have to
  keep track of two objects: a twisting\--invariant for the
  focus\--focus singularities (which is the same as the original one)
  and an invariant which keeps track of the normalization around the
  hyperbolic singularities. These objects should take into account the
  foliation type invariant replacing the polygon invariant in (i).

 \begin{figure}[h]
   \centering
   \includegraphics[width=0.35\linewidth]{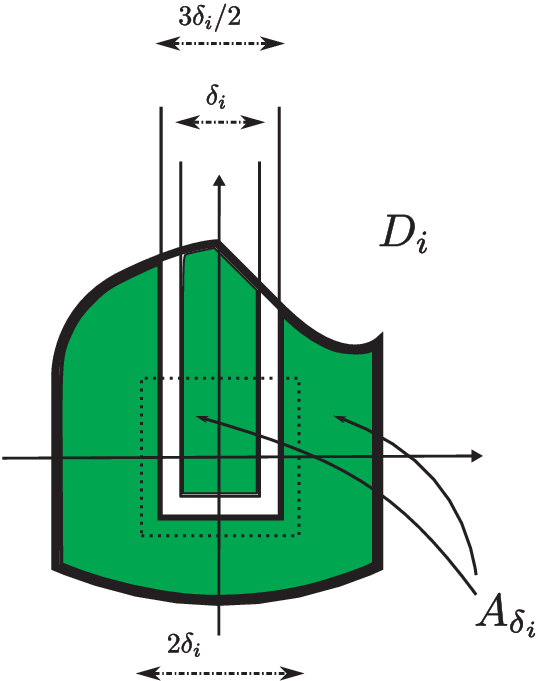}
   \caption{Symplectic gluing over polygon preimage. The point where
     the axes cross is a focus\--focus value. $D_i$ is a simply
     connected set containing this point.  The outer part of
     $A_{\delta_i}$ is a subset of $\mathbb{R}^2$ containing only
     regular values and on which we have a toric momentum map. The
     u-shaped region separates this region from a region of
     $\mathbb{R}^2$ which contains a focus\--focus value over which
     the system is not toric. These regions must be glued using a
     careful control of the smoothness and symplectic geometry of the
     system around the u-shaped region. Symplectic gluing arguments of
     this nature lie at the heart of the classification proofs
     (Theorem \ref{inventiones} and Theorem \ref{acta}).}
   \label{fig:spectra}
 \end{figure}

\item[(iii)] {\bf \emph{Symplectic gluing}}.  Once the invariants are
  defined, the symplectic gluing constructions outlined in the proof
  strategies of Theorems \ref{inventiones} and \ref{acta} need to be
  adapted. Indeed, the proofs need to incorporate the fact that in
  Step 1 we no longer have an invariant as simple as a polygon and
  hence the ``convenient covering'' of the polygon mentioned therein
  (\cite{PeVN2009,PeVN2011}) has to be changed. This convenient
  covering was such that it isolated the focus\--focus singularities
  of the integrable system, so that the preimage of the system over
  such element of the covering contained at most a focus\--focus
  singularity, which made the gluing construction feasible.

  Now we are presented with the fact that the image of the hyperbolic
  singularities in $F(M)$ can be a curve and the preimage contains
  many non\--isolated singularities.  How does one choose the elements
  $U_{\alpha}$ covering $\Delta$, so that the preimages
  $F^{-1}(U_{\alpha})$ can be glued symplectically and smoothly? In
  the case of semitoric systems outlined in the theorems, we were
  essentially gluing along regular points, but it is not clear if we
  can do this now: it seems quite likely that one needs to glue pieces
  taking care of the gluing along curves of singularities which, of
  course, must match.  One could probably use a covering of the curve
  of hyperbolic critical values, but we are presented with two
  problems:
  \begin{itemize}
  \item[(1)] \emph{Curves of the systems}: we do not know whether one
    can assume that for two given integrable systems, the
    aforementioned curves corresponding to each system, coincide. This
    was essential in the proof for semitoric systems. Moreover, does
    this information need to be included as an invariant?  Solving
    this problem will involve studying the singularities of the
    integral affine structure of the system.
  \item[(2)] \emph{Semiglobal Normal Forms}: we do not know how to
    glue because we do not have semiglobal normal forms around the
    singularities, these need to be proven first.
  \end{itemize}
  This poses a real challenge.

\end{itemize}

\subsection{\textcolor{black}{Non-proper systems}} \label{sec:proper}

The authors are working together with T.S. Ratiu on this part; the
fiber connectivity is proved in \cite{PeRaVN2011} and the construction
of the polygon invariant will be given in an upcoming article.

\subsubsection{Assumptions}  Suppose that $F \colon M \to
\mathbb{R}^2$ is an integrable system which fails to be semitoric only
because $J \colon M \to \mathbb{R}$ is not a proper map, but $F$ is a
proper map.  Concretely, this means that $(M,\, \omega)$ is a
connected symplectic $4$\--manifold equipped with two smooth functions
$J \colon M \to \R$ and $H \colon M \to \R$ such that:
\begin{itemize}
\item[(a)] $J$ is constant along the flow of the Hamiltonian vector
  field $\mathcal{H}_H$ generated by $H$ or, equivalently,
  $\{J,\,H\}=0$;
\item[(b)] for almost all points $p \in M$, the vectors
  $\mathcal{H}_J(p)$ and $\mathcal{H}_H(p)$ are linearly independent;
\item[(c)] $J$ generates a $2\pi$\--periodic flow, i.e., $J$ is the
  momentum map of an $S^1$\--action on $M$;
\item[(d)] $F$ is a proper map (but $J$ may not be a proper map);
\item[(e)] $F$ has only non\--degenerate and non-hyperbolic
  singularities.
\end{itemize}

The following is the classification problem of the section.

\begin{classproblem} \label{problemv2} \normalfont The problem has
  three parts.
  \begin{itemize}
  \item[{\rm (I)}] Give explicit constructions of ``analogues" of the
    symplectic invariants in Theorem \ref{inventiones} (which refer to
    semitoric systems only), for integrable systems satisfying
    assumptions (a)\--(e).
  \item[{\rm (II)}] Define an abstract list of all such possible
    invariants which occur in (I).
  \item[{\rm (III)}] Extend Theorems \ref{inventiones}, \ref{acta} to
    integrable systems satisfying (a)\--(e), using (I) and (II).
  \end{itemize}
\end{classproblem}

\subsubsection{Expected difficulties}  The difficulties we expect are as follows:

\begin{itemize}
\item[(i)] {\bf \emph{Fiber connectivity and the polygon invariant}}.
  In contrast with the standard semitoric setting, the fibers of $F$
  are not necessarily connected (see our paper \cite{PeRaVN2011} for
  the known results in this direction).  In order to guarantee that
  the fibers of $F$ are connected, in general, one can assume a
  non\--vertical tangency condition on the bifurcation diagram of the
  system, up to smooth deformations (see \cite[Theorem 1 and Theorem
  2]{PeRaVN2011}). This sufficient condition seems to us to be very
  close to necessary, but we do not know at this time if it can be
  weakened (there are many integrable systems with disconnected fibers
  that nevertheless do not satisfy this vertical tangency condition).

  Now, assuming that the fibers of $F$ are connected, it should be
  possible to construct a generalization of the polygon invariant,
  which will likely not be as simple as a family of polygons, but it
  will still be a rigid object that can be drawn in $\mathbb{R}^2$. We
  have investigated this case, and although we have no proofs yet, we
  think that analogue of the polygon invariant will be given in this
  case by a collection of regions of $\mathbb{R}^2$ (under some
  discrete group action), where each region is bounded by the graphs,
  epigraphs, or hypographs of two \emph{piecewise linear functions}.


  Note that, before constructing the polygon invariant, it is
  necessary to have a complete understanding of the image $F(M)$ of
  the system. It is precisely because we can now describe this image
  concretely (\cite[Theorem 3 and Theorem 4]{PeRaVN2011}) as bounded
  by the hypograph and epigraph of two lower/upper semicontinuous
  functions, that we can guess that the polygon invariant is going to
  be bounded by the graphs, epigraphs, or hypographs of two
  \emph{piecewise linear functions} (it should not be too difficult to
  do the transition from $F(M)$ to the polygon invariant, the proof
  should be very close to the original construction for semitoric
  systems).

  If, on the other hand, the map $F$ is proper but does not have
  connected fibers, then we do not know whether it is possible to
  construct a reasonable analogue of the polygon invariant. The
  original construction of the polygon invariant uses in an essential
  way the connectivity of the fibers of $F$. Even in the construction
  of the polytope corresponding to a toric integrable system $F$
  (i.e., $F$ is the momentum map of a Hamiltonian $2$\--torus action),
  the connectivity of the fibers of $F$ is intimately connected to the
  structure of $F(M)$ itself as a convex polytope. This is the context
  of the famous convexity theorem by Atiyah and Guillemin\--Sternberg
  (see Theorem \ref{theo:ags}), which says that $F(M)$ is, moreover,
  equal to the convex hull of the images of the fixed points of the
  action.

\item[(ii)] {\bf \emph{Other invariants}}. Until the analogue of the
  polygon invariant in (i) is constructed we do not see a way to
  anticipate what the other invariant which depends on the polygon
  (twisting\--index invariant) will be.  The cardinality, singular
  foliation type, and volume invariant should remain the same as for
  semitoric systems. In the case that $F$ has connected fibers we
  expect that the construction of these other invariants should not be
  too different from the case of semitoric systems.

\item[(iii)] {\bf \emph{Symplectic gluing}}.  The difficulties here
  will depend on the construction of the polygon. As in subsection
  \ref{sec:hyperbolic}, the constructions outlined in the proof
  strategies of Theorem \ref{inventiones} and Theorem \ref{acta} need
  to be adapted to deal with the fact that in Step 1 we no longer have
  an invariant as simple as a polygon, but rather a more complicated,
  but still ``rigid'' invariant which is close to a family of polygons
  inside of $\mathbb{R}^2$. The method of finding a suitable covering
  of a representative polygon of this family, as outlined in Step 1,
  should work. Of course, the covering here will be more complicated,
  but the fact that there are no hyperbolic singularities should make
  this step technically challenging but within reach.
\end{itemize}

\subsection{\textcolor{black}{Higher dimensional
    systems}} \label{sec:higher}

\subsubsection{Assumptions} When $M$ is a $6$\--dimensional symplectic
manifold, we think that the natural definition of a semitoric system
should be the following.

\begin{definition}
  A \emph{semitoric system} on a six\--dimensional symplectic manifold
  $M$ is a map $F \colon M \to \mathbb{R}^3$, $F=(J,\,H,\,K)$, where
  $J \colon M \to \mathbb{R}$, $H \colon M \to \mathbb{R}$, $K \colon
  M \to \mathbb{R}$ are smooth functions and the following conditions
  are satisfied.
  \begin{itemize}
  \item[(a)] $J,\,H,\,K$ are constant along the flows of the
    Hamiltonian vector fields they generate, i.e.,
$$
\{J,\,H\}=\{J,\,K\}=\{K,\,H\}=0;
$$
\item[(b)] for almost all points $p \in M$, the vectors
  $\mathcal{H}_J(p)$, $\mathcal{H}_H(p)$, $\mathcal{H}_K(p)$ are
  linearly independent;
\item[(c)] $J$ and $H$ generate $2\pi$\--periodic flows, i.e.,
  $(J,\,H)$ is the momentum map of an $(S^1)^2$\--action on $M$;
\item[(d)] $J$ and $H$ are proper maps;
\item[(e)] $F$ has only non\--degenerate singularities and
  non\--hyperbolic singularities.
\end{itemize}
\end{definition}

The following is the classification problem of the section.

\begin{classproblem} \label{problemv3} \normalfont The problem has
  three parts.
  \begin{itemize}
  \item[{\rm (I)}] Give explicit constructions of ``analogues" of the
    symplectic invariants in Theorem \ref{inventiones} (which refer to
    semitoric systems only), for integrable systems satisfying
    assumptions (a)\--(e).
  \item[{\rm (II)}] Define an abstract list of all such possible
    invariants which occur in (I).
  \item[{\rm (III)}] Extend Theorems \ref{inventiones}, \ref{acta} to
    integrable systems satisfying (a)\--(e), using (I) and (II).
  \end{itemize}
\end{classproblem}

\subsubsection{Expected difficulties.} The difficulties we expect are
as follows.

\begin{itemize}
\item[(i)] {\bf \emph{Fiber connectivity and the polygon
      invariant}}. A major difficulty is to define a $3$\--dimensional
  polytope invariant, because the singularities of $F$ are now more
  complicated. One can have singularities whose Williamson type
  (Definition \ref{def:zung}) is 
  $$(k_{\op{e}},\, k_{\op{h}},\,
  k_{\op{f}})=(1,\,0,\,1)\,\, \textup{or}\,\, (0,\,0,\,1), 
  $$
  for example. The
  singularities with the Williamson type $(0,\,0,\,1)$ have a
  focus\--focus component, but are no longer isolated as was the case
  for semitoric systems. This is already a major difference with the
  semitoric case in dimension four, and we believe it is a very
  substantial problem.
 
  The fundamental issues here are:
  \begin{itemize}
  \item \emph{what} is the stratification of the singular affine
    structure of the image $F(M)$, and
  \item can one describe concretely the singular stratum $\mathcal{S}$
    corresponding to the focus\--focus singularities?
  \end{itemize}
  The structure of $\mathcal{S}$ ought to be essential if one wants to
  cut $F(M)$ and, from it, construct a $3$\--dimensional polytope. At
  the moment, this is an outstanding and, we believe, a difficult
  problem.  In~\cite{Wa2012}, Wacheux proved that $\mathcal{S}$ does
  not contain any embedded circle~: it has to connect the stratum of
  transversally elliptic singularities. This opens the way to the
  construction of the polytope invariant.

  Note that in the case when the system $F \colon M \to \mathbb{R}^3$
  is the momentum map for a Hamiltonian action of a $3$\--torus, the
  polygon invariant is precisely the image $F(M)$, which was shown to
  be equal to the convex hull of the images of the fixed points of the
  action by the Atiyah-Guillemin\--Sternberg theorem (see Theorem
  \ref{theo:ags} and Figure \ref{fig:polytopes}).  In this case, the
  image of $F(M)$ is not ``cut''.

\item[(ii)] {\bf \emph{Other invariants}}. Once the polytope invariant
  is constructed, we expect difficulties, both conceptual and
  technical, in defining the other invariants with the possible
  exception of the singular foliation type invariant (which should be
  close to a twisted direct sum of the singular foliation type
  invariant in the semitoric case, with the additional component(s)
  corresponding to the other singularities).

  To arrive at a classification, both in the six\--dimensional, as
  well as the higher dimensional cases, it is necessary to understand
  the semiglobal classification of singularities. For instance, what
  are the semiglobal symplectic invariants near a singularity that has
  a focus\--focus and an elliptic component? This should not be too
  difficult, but as far as we know, it has not been done. See
  \cite[Section 5.2.3]{PeVN2012} for a summary and references of a
  semiglobal analysis of focus\--focus singularities in the
  $4$\--dimensional case, which was previously carried out by the
  second author.

  The twisting index invariant in this case should compare
  normalizations of the system near singularities which have
  focus\--focus components. An expected difficulty will come from the
  fact that in dimension 6 a focus-focus component is always coupled
  with either an elliptic or a regular component.  Thus these
  singularities come in families, which makes it difficult to isolate
  them or order them in some meaningful way (in the case of semitoric
  systems, the focus\--focus singularities were ordered according to
  the value of the $J$-component). Altogether, how to handle these
  issues is going to depend on how the polytope invariant is
  constructed.

\item[(iii)] {\bf \emph{Dimensions greater than $6$}}.  In part a)
  above, we only considered integrable systems on $6$\--dimensional
  manifolds.  After this case is dealt with, the analogous problem in
  any dimension would be to consider an integrable system given by
  smooth functions
$$
f_1,\, \ldots,\, f_{n-1}, \, f_n \colon M \to \mathbb{R}
$$
where each $f_1, \ldots, f_{n-1}$ generates a vector field with a
$2\pi$\--periodic flow. We shall call these systems \emph{semitoric
  systems with $n$ degrees of freedom} (what we have been calling
``semitoric systems'' are therefore the same as ``semitoric systems
with $2$ degrees of freedom").

\item[(iv)] {\bf \emph{Symplectic gluing}}.  A discussion analogous to
  that in subsections \ref{sec:hyperbolic} and \ref{sec:proper}
  applies. The construction of the covering of the $3$\--dimensional
  polytope by sets $U_{\alpha}$ will depend heavily on how the
  polytope is constructed and what the structure of the set of
  singularities inside of the polytope is.
\end{itemize}

\subsection{\textcolor{black}{Non-periodic systems}}
\label{sec:antiperiodic}

\subsubsection{Assumptions.} Suppose that $F \colon M \to \mathbb{R}^2$
is an integrable system which fails to be semitoric only because $J
\colon M \to \mathbb{R}$ is not a momentum map for a Hamiltonian
$S^1$\--action in $M$ (in other words, the Hamiltonian vector field
generated by $J$ does not have a $2\pi$\--periodic flow).  Concretely,
this means that $(M,\, \omega)$ is a connected symplectic
$4$\--manifold equipped with two smooth functions $J \colon M \to \R$
and $H \colon M \to \R$ such that:
\begin{itemize}
\item[(a)] $J$ is constant along the flow of the Hamiltonian vector
  field $\mathcal{H}_H$ generated by $H$ or, equivalently,
  $\{J,\,H\}=0$;
\item[(b)] for almost all points $p \in M$, the vectors
  $\mathcal{H}_J(p)$ and $\mathcal{H}_H(p)$ are linearly independent;
\item[(c)] $J$ may or may not generate a $2\pi$\--periodic flow (i.e.,
  $J$ may or may not be the momentum map of a Hamiltonian
  $S^1$\--action on $M$);
\item[(d)] $J$ is a proper map;
\item[(e)] $F$ has only non\--degenerate and non-hyperbolic
  singularities.
\end{itemize}

The following is the classification problem of the section.

\begin{classproblem} \label{problemv4} \normalfont The problem has
  three parts.
  \begin{itemize}
  \item[{\rm (I)}] Give explicit constructions of ``analogues" of the
    symplectic invariants in Theorem \ref{inventiones} (which refer to
    semitoric systems only), for integrable systems satisfying
    assumptions (a)\--(e).
  \item[{\rm (II)}] Define an abstract list of all such possible
    invariants which occur in (I).
  \item[{\rm (III)}] Extend Theorems \ref{inventiones}, \ref{acta} to
    integrable systems satisfying (a)\--(e), using (I) and (II).
  \end{itemize}
\end{classproblem}

\subsubsection{Expected difficulties.}

\begin{itemize}
\item[(i)] {\bf \emph{Fiber connectivity and the polygon invariant}}.
  Similarly to subsection \ref{sec:hyperbolic}, one should be able to
  define, from the image $F(M)$ of the integrable system $F$, a
  geometric invariant which resembles the polytope invariant. The
  fibers of $F$ may not be connected, which is essential in the
  construction of the polygon invariant for semitoric systems, as we
  have mentioned already. However, we know that these fibers are
  connected, for instance when the fibers of $J$ or the fibers of $H$
  are connected (\cite[Theorem 3.7]{PeRaVN2011}). The fibers of $F$
  can also be guaranteed to be connected by assuming that the
  bifurcation diagram of $F$ has no vertical tangencies, up to smooth
  deformations.

\item[(ii)] {\bf \emph{Polygon invariant, lack of momentum map, and
      other invariants}}. The fiber connectivity for $F$ is only a
  requirement in the construction of the polygon invariant; it allows
  us to describe $F(M)$ precisely, in terms of hypographs/epigraphs of
  lower/upper semicontinuous functions. However, in the construction
  of the polygon invariant it is essential that $J$ is the momentum
  map for a Hamiltonian $S^1$\--action on $M$. Indeed, this gives a
  ``vertical invariance'' to the image of $F(M)$, which allows one to
  cut $F(M)$ along vertical lines going through the focus\--focus
  values of $F$, and hence to construct a polygon from $F(M)$. If we
  drop this periodicity assumption on $J$, the construction of the
  polygon needs to be rethought; we believe this to be a challenging
  task (at least a priori). The work of Leung and Symington \cite{ls}
  is directly relevant to these problems, in case the fibers are
  connected.

\item[(iii)] {\bf \emph{Symplectic gluing}}.  A discussion analogous
  to that of subsections \ref{sec:hyperbolic} and \ref{sec:proper}
  applies.

\end{itemize}

\subsection{\textcolor{black}{Degenerate
    systems}} \label{sec:degenerate}

This part of the program, as was the case with subsection
\ref{sec:hyperbolic}, is wide open.  We do not know how realistic it
is to achieve a classification for systems that have degenerate
singularities.  Degenerate singularities appear in many systems, so
the study of degenerate systems should be considered an important
case. Some important systems, such as the three wave interaction (see
subsection \ref{sec:first}), are known to have degenerate
singularities.

We believe that the work of Garay~\cite{garay,garay2}, Sevenheck and
van Straten~\cite{sevenheck,sevenheck2} on Lagrangian rigidity, and
Zung's work on analytic Birkhoff normal forms \cite{Zu2005}, make very
substantial advances towards the understanding of degenerate
singularities.

The one\--dimensional $\Cinf$ smooth case was done by de
Verdi{\`e}re. It involves Gauss-Manin and techniques by Malgrange
related to Milnor fibers, and usual singularity theory.

\subsubsection{Assumptions.}  Suppose that $F \colon M \to
\mathbb{R}^2$ is an integrable system which fails to be semitoric only
because it may have degenerate singularities.  Concretely, this means
that $(M,\, \omega)$ is a connected symplectic $4$\--manifold equipped
with two smooth functions $J \colon M \to \R$ and $H \colon M \to \R$
such that:
\begin{itemize}
\item[(a)] $J$ is constant along the flow of the Hamiltonian vector
  field $\mathcal{H}_H$ generated by $H$ or, equivalently,
  $\{J,\,H\}=0$;
\item[(b)] for almost all points $p \in M$, the vectors
  $\mathcal{H}_J(p)$ and $\mathcal{H}_H(p)$ are linearly independent;
\item[(c)] $J$ generates a $2\pi$\--periodic flow;
\item[(d)] $J$ is a proper map;
\end{itemize}

Note that unlike in previous subsections, there is no item (e)
above. If there are degenerate singularities there is no a priori
reason to single out hyperbolic singularities. We have the following
outstanding preliminary problem, which is self\--contained.

\begin{coreproblem} \label{ppp} \normalfont Determine what kinds of
  singularities can be allowed for a system satisfying (a)\--(d), so
  that one can assure that the system has connected fibers.
\end{coreproblem}

The following is the classification problem of the section.

\begin{classproblem} \label{problemv5} \normalfont The problem has
  three parts.
  \begin{itemize}
  \item[{\rm (I)}] Give explicit constructions of ``analogues" of the
    symplectic invariants in Theorem \ref{inventiones} (which refer to
    semitoric systems only), for integrable systems satisfying
    assumptions (a)\--(d), and which only have the types of
    singularities in Problem \ref{ppp}.
  \item[{\rm (II)}] Define an abstract list of all such possible
    invariants which occur in (I).
  \item[{\rm (III)}] Extend Theorems \ref{inventiones}, \ref{acta} to
    integrable systems satisfying (a)\--(d), using (I) and (II).
  \end{itemize}
\end{classproblem}

\subsubsection{Expected difficulties.} Very few results are known for degenerate singularities (see
\cite[Section 4.2.2]{PeVN2012} for a brief discussion and
references). Eliasson's linearization theorem, on which much of the
modern theory of finite dimensional completely integrable Hamiltonian
systems relies (and in particular the results of \cite{PeVN2009,
  PeVN2011, PeVN2010, PeRaVN2011}), holds only for
non\--degenerate singularities.

However, degenerate singularities occur often in systems arising from
physics, so it is important to develop a theory which covers systems
with degenerate singularities. A first step is to investigate whether
any part of Eliasson's theorem can be rescued in the degenerate case
(and the hope is small of getting a result which is not much weaker
than Eliasson's original theorem).

There has been, however, some work for degenerate singularities in
\cite[Section 4.2.2]{PeVN2012}, and it is conceivable that the
classification theorems we proved in \cite{PeVN2009, PeVN2011} could
be extended to cover the degenerate case by reformulating these
results. This is definitely a difficult, but very important, part of
the classification program announced in this paper. This part stands
somewhat independently of the rest of the program and is the one for
which we have fewer ideas or expectations of where it can lead.

We do not know how having degenerate singularities will affect the
construction of the symplectic invariants. First one needs to classify
some subclass of degenerate singularities
(cf. Arnold\--Gusein-Zade\--Varchenko's books
\cite{ArGuVa1985})

We expect this part of the program to be of great interest for the
physics community.

\subsection{\textcolor{black}{Self-contained topological and geometric
    questions}} \label{sec:topology}

There are many self-contained questions in this program that are of
interest on their own right. The following gives a small sample. In
addition, answering some of them is a prerequisite, as far as we can
tell, for the construction of the polygon invariant in the previous
sections. This is the case for Problem \ref{pr1a} and Problem
\ref{pr1b} below; therefore it can be a good strategy to start by
solving these problems in the context of the corresponding section
(i.e., hyperbolic systems, non-proper systems, non-periodic systems
etc.).

Because the Fomenko school has done work on topological aspects of
integrable systems, we would like to point out how the following
questions, though topological, differ from those they answered.  The
Fomenko school gave an extensive treatment of the \textit{topological}
properties of integrable systems viewed as fibrations over
$\mathbb{R}^n$ (we refer to \cite{PeRaVN2011} for references); this
work has exerted a major influence on several parts of modern
mathematics, including our own work.

\begin{figure}[h]
  \centering
  \includegraphics[width=0.5\textwidth]{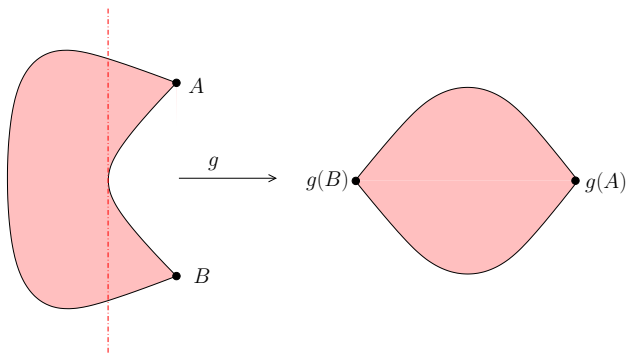}
  \caption{Suppose that the bifurcation set $\Sigma_F$ of $F$ consists
    precisely of the boundary points in the left figure (which depicts
    $F(M)$). The diffeomorphism $g$ transforms $F(M)$ to the region on
    the right hand side of the figure, in order to remove the original
    vertical tangencies on $\Sigma_F$.}
  \label{fig:diffeo}
\end{figure}

The goal and results of the \emph{symplectic} bifurcation theory
introduced in \cite{PeRaVN2011} are fundamentally different from those
of the Fomenko school.  For instance, from the Fomenko school
view-point, if the fibration given by an integrable system has
disconnected fibers, one may try to slightly modify the system, to
consider a topologically equivalent system with connected fibers
(``topologicaly equivalent" in the sense that sets of connected
components of the fibers of both systems coincide). In general, such a
change does not respect the symplectic structure, so the modified
system can be topologically equivalent to the original one but its
symplectic geometry and dynamical behavior can be different. Since we
are interested in quantization and spectral theory of the integrable
system, such an operation defeats our purpose: the symplectic
structure cannot be altered.

In fact, the fibrations of two integrable systems may be globally the
same at a topological level, while their symplectic invariants are
inequivalent \cite{PeVN2009}. From the point of view of symplectic and
spectral theory the two systems have little to do with each other:
there is no global isomorphism which preserves the integrable system
with the symplectic structure, which in turn changes the way
quantization and spectral theory are carried out.

The development of the symplectic bifurcation theory in
\cite{PeRaVN2011} requires the introduction of methods to construct
Morse-Bott functions which, from the point of view of symplectic
geometry, behave well near the singularities of integrable systems.
These methods use Eliasson's theorems on linearization of
non-degenerate singularities as well as the symplectic topology of
integrable systems, topics that have been developed by many.

\begin{theorem}[Pelayo-Ratiu-V\~ u Ng\d oc] \label{theo:main-connectivity0} Suppose that $(M, \omega)$ is
  a compact connected symplectic four-manifold. Let $F \colon M
  \to\mathbb{R}^2$ be a non-degenerate integrable system without
  hyperbolic singularities. Denote by $\Sigma_F$ the bifurcation set
  of $F$.  Assume that there exists a diffeomorphism $g \colon F(M)
  \to \mathbb{R}^2$ onto its image such that $g(\Sigma_F)$ does not
  have vertical tangencies (see Figure~\ref{fig:diffeo}).  Then $F$
  has connected fibers.
\end{theorem}

\begin{figure}[h]
  \centering
  \includegraphics[width=0.8\textwidth]{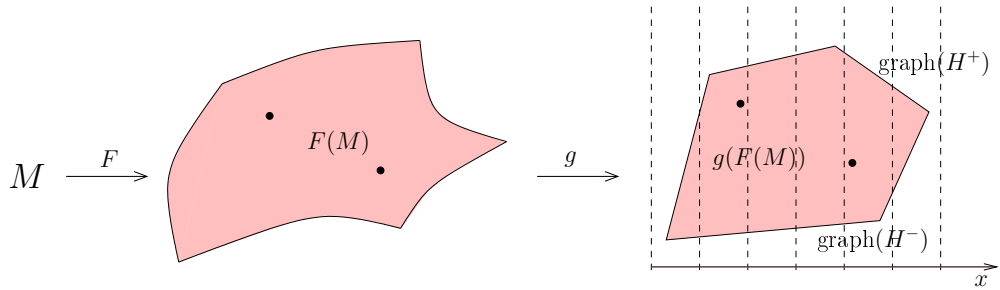}
  \caption{Description of the image of an integrable system. The image
    is first transformed to remove vertical tangencies, and then it
    can be described as a region bounded by two graphs.}
  \label{fig:fig1}
\end{figure}

\begin{coreproblem} \label{pr1a} \normalfont Prove Theorem
  \ref{theo:main-connectivity0} under the assumptions of each of the
  previous sections of the paper.  For instance, in higher dimensions:
  suppose that $(M, \omega)$ is a compact connected symplectic
  $2n$-manifold. Let $F \colon M \to\mathbb{R}^n$ be a non-degenerate
  integrable system without hyperbolic singularities. Denote by
  $\Sigma_F$ the bifurcation set of $F$.  Assume that there exists a
  diffeomorphism $g \colon F(M) \to \mathbb{R}^n$ onto its image such
  that $g(\Sigma_F)$ does not have ``vertical tangencies''.  Does $F$
  always have connected fibers?  The first problem is to determine
  precisely what we mean by ``vertical tangency'' for a domain in
  $\mathbb{R}^n$ when $n>2$.
\end{coreproblem}

\begin{coreproblem} \label{pr1b} \normalfont There is also an analogue
  of Theorem \ref{theo:main-connectivity0} for non-compact manifolds
  in \cite{PeRaVN2011}, so we can also pose the problem above for
  non-compact manifolds.
\end{coreproblem}

\begin{problem}
  \normalfont Find formulas, or recipes, to compute the symplectic
  invariants constructed in the previous section. For instance, is
  there a general strategy to compute some terms of the Taylor series
  invariant (i.e., the singularity type invariant)? For instance, if
  the system is toric, then the image of the momentum map, which is
  the simplest case of a polygon invariant, can be computed in terms
  of the images of the fixed points of the action. Can one write a
  recipe to compute it in each of the cases described in the previous
  sections of this paper?
\end{problem}

\begin{problem}
  \normalfont Are there Duistermaat\--Heckman formulas suitable for
  the polygon invariant?
\end{problem}

The following result describes the image of an integrable system.

\begin{theorem}[Pelayo-Ratiu-V\~ u Ng\d oc]
  \label{theo:main-image0}
  Suppose that $(M, \omega)$ is a compact connected symplectic
  four-manifold. Let $F \colon M \to \mathbb{R}^2$ be a non-degenerate
  integrable system without hyperbolic singularities. Denote by
  $\Sigma_F$ the bifurcation set of $F$.  Assume that there exists a
  diffeomorphism $g \colon F(M) \to \mathbb{R}^2$ onto its image such
  that $g(\Sigma_F)$ does not have vertical tangencies (see
  Figure~\ref{fig:diffeo}). Then:
  \begin{itemize}
  \item[{\rm (1)}] the image $F(M)$ is contractible and the
    bifurcation set can be described as
  $$\Sigma_F=\partial  (F(M)) \sqcup \mathcal{F},$$ where $\mathcal{F}$ is
  a finite set of rank $0$ singularities which is contained in the
  interior of $F(M)$;
\item[{\rm (2)}] let $(J,\,H):=g\circ F$ and let $J(M)=[a,\,b]$. Then
  the functions $H^{+},\, H^{-} \colon [a,\,b] \to \R$ defined
  by $$H^{+}(x) := \max_{J^{-1}(x)} H$$ and $$H^{-}(x) :=
  \min_{J^{-1}(x)} H$$ are continuous and $F(M)$ can be described
  as $$F(M) = g^{-1}(\op{epi}(H^{-}) \cap \op{hyp}(H^+)).$$
\end{itemize}
\end{theorem}

  \begin{figure}[h]
    \centering
    \includegraphics[width=0.6\textwidth]{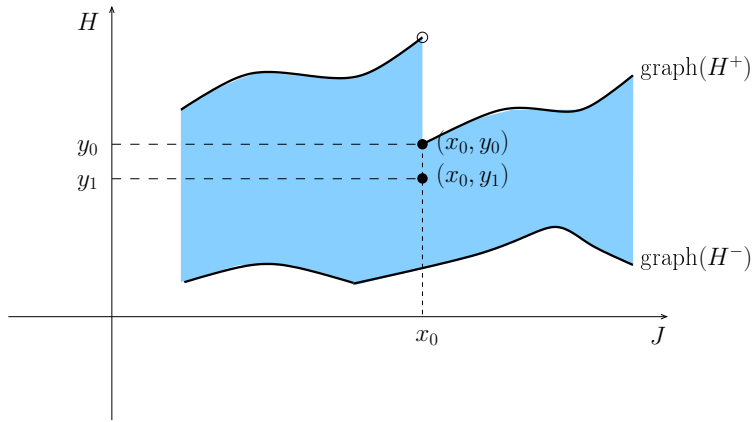}
    \caption{Image $F(M)$ of an integrable system with
      non\--degenerate and non\--hyperbolic singularities. The point
      $(x_0,\,y_1)$ is contained in the interior of $F(M)$, while the
      point $(x_0,\,y_0)$ is contained in the topological boundary of
      $F(M)$.}
    \label{fig:lowersemicontinuity}
  \end{figure}

  \begin{coreproblem} \label{pr1b2} \normalfont Prove Theorem
    \ref{theo:main-image0} under the assumptions of each of the
    previous sections of the paper.  For instance, in higher
    dimensions: suppose that $(M, \omega)$ is a compact connected
    symplectic $2n$-manifold. Let $F \colon M \to\mathbb{R}^n$ be a
    non-degenerate integrable system without hyperbolic
    singularities. Denote by $\Sigma_F$ the bifurcation set of $F$.
    Assume that there exists a diffeomorphism $g \colon F(M) \to
    \mathbb{R}^n$ onto its image such that $g(\Sigma_F)$ does not have
    ``vertical tangencies''.  Does $F$ always have connected fibers?
  \end{coreproblem}

  \begin{coreproblem} \label{pr1c} \normalfont There is also an
    analogue of Theorem \ref{theo:main-image0} for non-compact
    manifolds in \cite{PeRaVN2011}, so we can also pose the problem
    above for non-compact manifolds.
  \end{coreproblem}

\begin{figure}[h]
  \centering
  \includegraphics[width=0.75\textwidth]{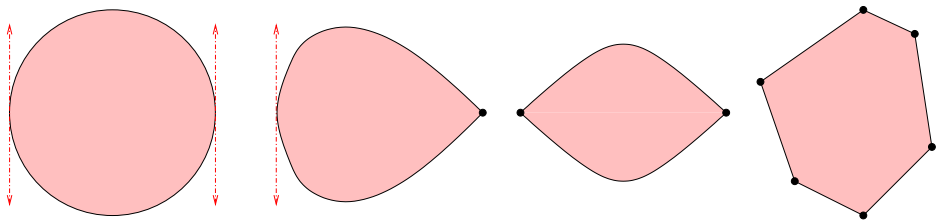}
  \caption{A disk, a disk with a conic point, a disk with two conic
    points, or a compact convex polygon.}
  \label{fig:4images}
\end{figure}

As a matter of fact, some cases in which vertical tangencies are
present, can also be addressed. In the following result a neighborhood
of a \emph{conic point} is, by definition, locally diffeomorphic to
some proper cone $C_{\alpha,\beta}$.

\begin{theorem}[Pelayo\--Ratiu\--V\~ u Ng\d oc]
  \label{theo:striking}
  Suppose that $(M, \omega)$ is a compact connected symplectic
  four-manifold. Let $F \colon M \to \mathbb{R}^2$ be a non-degenerate
  integrable system without hyperbolic singularities.  Assume that
  \begin{enumerate}
  \item[{\rm (a)}] the interior of $F(M)$ contains a finite number of
    critical values;
  \item[{\rm (b)}] there exists a diffeomorphism $g$ such that
    $g(F(M))$ is either a disk, a disk with a conic point, a disk with
    two conic points, or a compact convex polygon (see Figure
    \ref{fig:4images}).
  \end{enumerate}
  Then the fibers of $F$ are connected.
\end{theorem}

\begin{problem}
  \normalfont Prove an analogue of Theorem \ref{theo:striking} in
  dimensions $2n\ge 6$.
\end{problem}

We next pose three, in our opinion, important moduli problems which
are independent of the classification program.

\begin{problem} \label{pr1} \normalfont Describe the topology of the
  moduli space of \emph{symplectic toric systems} (using Theorem
  \ref{delzantseminal}).  The moduli of semitoric systems of a given
  dimension may be viewed inside of the set all of integrable systems
  in that dimension -- can one tell how ``big'' it is?
\end{problem}

As far as we know there are no definite results regarding Problem
\ref{pr1}. There is collection of notes about this problem by
A.R. Pires and the first author with several preliminary observations.
The following should be a much more challenging problem.

\begin{problem} \normalfont Describe the topology of the
  moduli space of semitoric systems (using Theorem \ref{inventiones}
  and Theorem \ref{acta}).
\end{problem}

The following problem is vaguely defined, but probably worth
exploring.

\begin{problem} \label{pr2} \normalfont Can one define a suitable
  notion of ``limit'', such that an integrable system with connected
  fibers in dimension four could be obtained as a ``limit" of a
  sequence of semitoric systems (at least for some type of
  non\--semitoric integrable systems)?
\end{problem}

The following problem has the same flavor as Problem \ref{pr2}, and is
motivated by an informal, famous question of Anatole Katok \cite{Ka}
in low dimensions (meaning $2$ for maps, $3$ for flows): ``Is every
conservative dynamical system that has zero topological entropy a
limit of integrable systems?"

\begin{problem}
  \normalfont Consider the set of integrable systems:
$$
\mathcal{I}:=\{ F \colon M \to \mathbb{R}^2 \,\,\, | \,\,\, F\,\,
\textup{is an integrable system}\} \subset \op{C}^{\infty}(M,\,
\mathbb{R}^2).
$$
Similarly one can define $\mathcal{I}_{\textup{toric}}$ and
$\mathcal{I}_{\textup{semitoric}}$, by requiring, respectively, the
integrable systems in $\mathcal{I}$ to be toric or semitoric.
\begin{itemize}
\item[(a)] What is the closure of $\mathcal{I} \subset
  \op{C}^{\infty}(M,\, \mathbb{R}^2)$, in some adequate topology to be
  defined?
\item[(b)] One can ask the same question for
  $\mathcal{I}_{\textup{toric}}$ and
  $\mathcal{I}_{\textup{semitoric}}$.
\end{itemize}
The same problem can be formulated for the systems in Sections
\ref{sec:hyperbolic}, \ref{sec:proper}, \ref{sec:higher},
\ref{sec:antiperiodic}, and \ref{sec:degenerate}.
\end{problem}

We conclude by raising some broad problems about semitoric systems.

\begin{problem} \label{pppp} \normalfont Study the homotopy groups,
  homology, equivariant cohomology of symplectic manifolds admitting
  semitoric systems. The computations of these groups should be in
  terms of the list of invariants of the system (see Theorems
  \ref{inventiones} and \ref{acta}). What restriction does the
  existence of a semitoric system pose on the topology of the
  manifold?
\end{problem}

\begin{problem} \label{pppp2} \normalfont Let $M$ be a compact
  symplectic manifold equipped with an effective Hamiltonian
  $S^1$\--action with momentum map $J \colon M \to \mathbb{R}$.  Under
  which conditions can one find a smooth function $H \colon M \to
  \mathbb{R}$ such that $F:=(J,\,H) \colon M \to \mathbb{R}^2$ is a
  semitoric system as in Section \ref{sec:semitoric}?
\end{problem}

\subsection{\textcolor{black}{Collaborative
    work}} \label{sec:collaborative}

There are three levels we suggest to go from the material described in
the previous sections to a symplectic classification of completely
integrable Hamiltonian systems.
\\
\\
{\bf Level 1} (\emph{Putting the classifications together}). This
level should involve no major technical or conceptual difficulties. We
see it as a serious book-keeping and organizational challenge: the
effort will be into coordinating the results obtained by the authors
that have respectively worked in subsections \ref{sec:proper},
\ref{sec:higher}, \ref{sec:antiperiodic}, \ref{sec:degenerate},
\ref{sec:topology}, and putting their classification results together.
As we mentioned in the introduction, we believe that, for the most
part, each of these sections deals with difficulties that are somewhat
unrelated. Precisely because of the likely ``independence", the
general answer in \ref{sec:collaborative} should be a ``direct sum''
of the answers provided in the subsections \ref{sec:hyperbolic},
\ref{sec:proper}, and \ref{sec:antiperiodic}.

At this level, putting the classifications obtained in subsections
\ref{sec:hyperbolic}, \ref{sec:proper}, \ref{sec:degenerate} together,
we expect to arrive at a result of the following type:
\\
\\
{\bf \emph{Goal: Symplectic Classification Theorem for Semitoric
    Systems}}.  An integrable system
$$
F:=(f_1,\, \ldots,\, f_{n-1},\, f_n) \colon M \to \mathbb{R},
$$
where each $f_1,\, \ldots, f_{n-1}$ generates a vector field with a
$2\pi$\--periodic flow is characterized -- in the sense of uniqueness
and existence as in Theorems \ref{inventiones} and \ref{acta} -- up to
symplectic isomorphisms preserving the system, by a list of symplectic
invariants which are explicitly defined (in subsections
\ref{sec:hyperbolic}, \ref{sec:proper}, and \ref{sec:antiperiodic}).
\\
\\
{\bf Level 2} (\emph{Complexity two semitoric systems}). We call
\emph{complexity of an integrable} system the number of non\--periodic
components (so a semitoric system has complexity one). This level
holds the major conceptual and technical jump of the program and will
rely on techniques analogous to those developed in subsection
\ref{sec:higher}. We believe that going from the four to the
six\--dimensional case already contains the essential difficulties and
those need to be clarified first. The type of theorem that we expect
to have is the following.
\\
\\
{\bf \emph{Goal: Symplectic Classification Theorem for Complexity Two
    Semitoric systems}}.  An integrable system
$$
F:=(f_1,\, \ldots,\, f_{n-2},\, f_{n-1}, \, f_n) \colon M \to \mathbb{R},
$$
where each $f_1,\, \ldots, f_{n-2}$ generates a vector field with a
$2\pi$\--periodic flow is characterized -- in the sense of uniqueness
and existence as in Theorems \ref{inventiones} and \ref{acta} -- up to
symplectic isomorphisms preserving the system, by a list of symplectic
invariants which are explicitly defined and analogous to the ones in
subsection \ref{sec:higher} (in terms of the invariants in subsections
\ref{sec:hyperbolic}, \ref{sec:proper}, and \ref{sec:antiperiodic}).
Related to this goal is, in the context of Hamiltonian torus actions,
the work of Karshon and Tolman on \emph{complexity one} Hamiltonian
torus actions. Their approach might be useful at this level.
\\
\\
{\bf Level 3} (\emph{Increasing the complexity of a semitoric system
  to arrive to a general integrable systems}).  Here there is another
really substantial step: to proceed inductively using the theorem in
level 2 above, to obtain a general classification theorem of the
following form.
\\
\\
{\bf \emph{Goal: Symplectic Classification Theorem for Integrable
    Systems}}.  An integrable system (or at least a large class of
integrable systems)
$$
F:=(f_1,\, \ldots,\, f_{n-1}, \, f_n) \colon M \to \mathbb{R},
$$
is characterized -- in the sense of uniqueness and existence as in
Theorems \ref{inventiones} and \ref{acta} -- up to symplectic
isomorphisms preserving the system, by a list of symplectic invariants
which are explicit (possibly defined recursively in terms of
invariants of lower complexity systems).
\\
\\
We believe that it is possible to achieve a classification in terms of
explicit invariants, but these are likely to be defined inductively in
terms of invariants of lower complexity systems.

Even if one deals with the quite rigid case of integrable systems
coming from torus actions, the goal seems still very challenging, and
at this point we cannot provide much information about it.  However,
once the the six\--dimensional case in Section~\ref{sec:higher} is
understood, one should be able to compare the invariants therein with
the invariants of semitoric systems in Theorems \ref{inventiones},
\ref{acta}, and conjecture an approach to go from six\--dimensional
manifolds to eight\--dimensional etc. In this sense, the following
problem is a (possible) preliminary step towards the goal above.

\begin{problem}
  \normalfont Describe to the greatest possible extent the symplectic
  invariants built in Section \ref{sec:higher} for semitoric systems
  on $6$\--dimensional manifolds, in terms of symplectic invariants of
  semitoric systems of $4$\--dimensional symplectic manifolds (in
  Theorems \ref{inventiones}, \ref{acta}).
\end{problem}

\section{Spectral theory of integrable
    systems} \label{sec:spectral} 
    Our original motivation to develop a symplectic theory of semitoric
systems in the papers above was that, at least from our point of view,
it is a prerequisite for studying the quantization and spectral theory
of integrable systems.

\subsection{\textcolor{black}{Inverse spectral theory}}

Let $M$ be a symplectic manifold of dimension $2n$ and
$\mathcal{H}_{\hbar}$ a family of Hilbert spaces, $\hbar>0$,
associated to $M$ in such a way that linear operators on
$\mathcal{H}_\hbar$ admit principal symbols in ${\rm C}^{\infty}(M)$.

For example, such a family of Hilbert spaces can be obtained by
geometric quantization with complex polarization (\cite{ChPeVN2012}),
or by hand like in quantum mechanics (\cite{PeVN2010}), or if $M$ is a
cotangent bundle minus the zero section some completion of the
Schwartz space.  The simplest example is $\op{L}^2(X)$ when $M={\rm
  T}^*X$, a case in which there is no dependance on $\hbar$.

 \begin{figure}[h]
   \centering
   \includegraphics[width=0.6\linewidth]{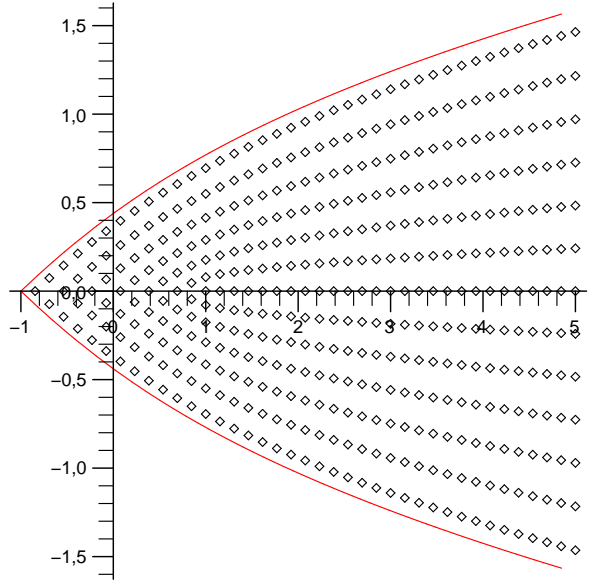}
   \caption{Spectrum of the quantum coupled spin\--oscillator.}
   \label{fig:spectra2}
 \end{figure}

\begin{definition}
\label{defi:quantum}
  A \emph{quantum integrable system} on $M$ consists of $n$
  semiclassical operators
$$
\hat{f}_1=(\hat{f}_{1,\hbar}) , \ldots , \hat{f}_n=(\hat{f}_{1,\hbar})
$$
acting on $\mathcal{H}_\hbar$ which commute (i.e.,
$[\hat{f}_{i,\hbar},\, \hat{f}_{j,\hbar}] = 0$ for all
$i,\,j=1,\ldots, n$, for all $\hbar$) and whose principal symbols
$f_1,\ldots, f_n$ form an integrable system on $M$.
\end{definition}

Instead of a single spectrum one considers the joint spectrum of
$\hat{f}_1 ,\ldots , \hat{f}_n$ (assumed discrete), i.e. the sequence
of points as $\hbar$ goes to $0$, given by:
$$
\op{Jspec}(\hat{f}_{1,\hbar},\ldots ,\hat{f}_{n,\hbar})=
\left\{(\lambda_{1,\hbar},\ldots ,\lambda_{n,\hbar})\in \mathbb{C}^n
  \,\,\left|\,\, \bigcap_{j=1}^n
    \op{ker}(\hat{f}_{j,\hbar}-\lambda_{j,\hbar}{\rm I}) \neq 0
  \right.\right\}.
$$

\begin{question}
  \normalfont Does the joint spectrum determine the integrable system
  $f_1,\ldots,f_n$ up to symplectic isomorphism?
\end{question}

\begin{example}(Quantum coupled spin\--oscillator) Recall the
  classical coupled spin\--oscillator model in Example
  \ref{ex:sp}. Now we are interested in its quantum version.  For any
  $\hbar>0$ such that $2=\hbar(n+1)$, for some non-negative integer
  $n\in\N$, let $\mathcal{H}$ denote the standard $n+1$-dimensional
  Hilbert space quantizing the sphere $S^2$ (see \cite{PeVN2010} for
  details). Consider on $\mathbb{R}^2$ coordinates $(u,\,v)$ and on
  $S^2$ coordinates $(x,\,y,\,z)$. The quantization of $x,\,y,\,z$ is
  given by restricting, respectively, the following operators to
  $\mathcal{H}$:
  \begin{eqnarray}
    \hat{x}:=\frac{\hbar}{2}(a_1a_2^*+a_2a_1^*), \qquad
    \hat{y}:=\frac{\hbar}{2\ii }(a_1a_2^*-a_2a_1^*), \qquad
    \hat{z}:=\frac{\hbar}{2}(a_1a_1^*-a_2a_2^*). \nonumber
  \end{eqnarray}
  where
$$a_i:=\frac{1}{\sqrt{2 \hbar}} \Big( \hbar
\frac{\partial}{\partial x_j} +x_j\Big),\,\,\,\,i=1,2.
 $$
 The unbounded operators
   $$\op{Id} \otimes \Big(
   -\frac{\hbar^2}{2} \frac{\op{d}^2}{\op{d}u^2} +\frac{u^2}{2} \Big)
   + (\hat{z} \otimes \op{Id})
   $$
   and
   $$
   \frac{1}{2}(\hat{x}\otimes u + \hat{y} \otimes
   (\frac{\hbar}{\ii}\frac{\partial}{\partial u}),
   $$
   which we call $\hat{J}$ and $\hat{H}$ respectively, on the Hilbert
   space $$\mathcal{H} \otimes \op{L}^2(\R)\subset \op{L}^2(\R^2)
   \otimes \op{L}^2(\R),$$ are self\--adjoint and commute. The
   spectrum of $\hat{J}$ is discrete and consists of eigenvalues in
   $\hbar(\frac{1-n}{2}+\N)$. This is natural way to quantize the
   coupled\--spin oscillator. The joint spectrum for a fixed value of
   $\hbar$ is depicted in Figure \ref{fig:spectra2}.
 \end{example}

\begin{figure}[h]
  \centering
  \includegraphics[width=0.5\textwidth]{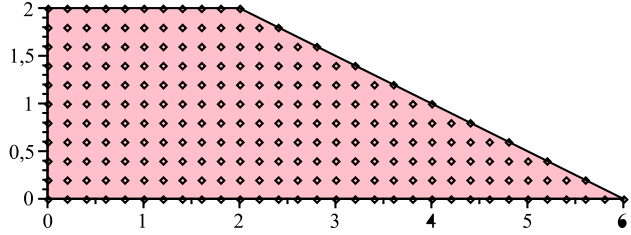}
  \caption{``Model image" of the spectrum of a normalized quantum
    toric integrable system. }
  \label{fig:normal}
\end{figure}

The question of whether one can recover a system from its spectrum is
fundamental in the theory of integrable systems.  The
\emph{symplectic} classification in terms of concrete invariants
described in the previous section serves as a tool to quantize
semitoric systems.

 \begin{figure}[h]
   \centering
   \includegraphics[width=0.5\textwidth]{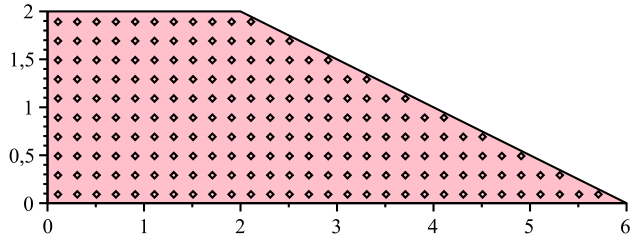}
   \caption{The figure shows the spectra in Figure \ref{fig:normal}
     with a so called ``metaplectic correction''.  Introducing a
     metaplectic correction refers to twisting the prequantum bundle
     (or its powers) by a half-form bundle. The metaplectic correction
     allows to obtain an easier control of the subprincipal terms in
     the semiclassical limit, see our article \cite{ChPeVN2012}.  }
   \label{fig:deformed_metaplectic}
 \end{figure}

 In Delzant's theory, the image of the momentum map, for a toric
 integrable action, completely determines the system. In the quantum
 theory, the image of the momentum map is replaced by the joint
 spectrum. Can one determine the system from the joint spectrum? In
 this vast program, one can ask the following question, less ambitious
 but still quite spectacular: does the semiclassical (i.e., for a
 sequence of small values of $\hbar$) joint spectrum determine the
 underlying classical system? This was one of the primary goals of the
 program described in this paper:
 \\
 \\
 {\bf \emph{Spectral Goal of Program for Integrable Systems}}:
 \emph{Show that from the semiclassical joint spectrum of a quantum
   integrable system one can recover the integrable system of
   principal symbols, up to symplectic isomorphisms} (a program
 proposal to arrive at this result is outlined in subsection
 \ref{strategy_quantum}).
 \\
 \\
 We believe that it is possible to make definite progress towards
 achieving this goal. In our article \cite{PeVN2012} we formulated
 this goal as a conjecture in the case of semitoric systems
 (\cite[Conjecture 9.1]{PeVN2012}).

 As it is evident from the strategy we describe below to achieve this
 goal, the achievement of the symplectic program described Section
 \ref{sec:symplectic} will to a large extent determine the success of
 this part of the program. At this stage, we can formulate a general
 conjecture.  Note that the notion of quantum system is not
 universally agreed upon and the conjecture may also hold knowing the
 semiclassical joint spectrum up to some order. Therefore, the
 following statement should be understood with a certain flexibility
 depending on the context. For instance, in the case of toric systems,
 a very complete answer can be given, which has as corollary a
 quantization theorem.

 \begin{conjecture}[Inverse Spectral Conjecture for Integrable
   Systems] \label{con} From the semiclassical joint spectrum of a
   quantum integrable system one can completely recover the integrable
   system given by the principal symbols, up to symplectic
   isomorphisms.
 \end{conjecture}

 We have proved this conjecture in the case of toric systems (Theorem
 \ref{theo:inverse-spectral}) in our recent paper \cite{ChPeVN2012}.
 In the case that $M$ is a $2$\--dimensional cotangent bundle and
 $\hat{f}_1$ is a pseudodifferential operator, the conjecture was
 proved by V\~{u} Ng\d{o}c~\cite{san-inverse}. The case when $M$ is a
 $2$\--dimensional symplectic manifold and $\hat{f}_1$ is a Toeplitz
 operator is a field of current interest. In this direction, a local
 quantum normal form near elliptic singularities has been obtained by
 Le Floch~\cite{lefloch}.  Note that in $2$\--dimensions, a
 singularity cannot have focus\--focus components.

 \subsection{\textcolor{black}{An approach to Conjecture \ref{con} for
     semitoric systems}}
 \label{strategy_quantum}

 We begin with the following preliminary remarks.

 \begin{itemize}
 \item[(i)] \emph{Toeplitz operators (on any manifold) versus
     pseudodifferential operators (on cotangent bundles)}. A first
   difficulty is to define precisely the space of quantum operators
   that are considered.  We believe that the natural setting to
   formulate this question is semiclassical Toeplitz
   operators. Nevertheless, the approach below is based on ideas
   coming from the theory of pseudodifferential operators.

   Pseudodifferential quantization implies that the classical phase
   space be a cotangent bundle $M={\rm T}^*X$, which is a serious
   restriction from the point of view of symplectic geometry.
   However, it is well known that a microlocal level the algebra of
   pseudodifferential operators is microlocally equivalent to the
   algebra of Toeplitz operators, and one can
   reasonably expect that the techniques that were developed for
   pseudodifferential operators have their analogues for Toeplitz
   operators.
 \item[(ii)] \emph{Semitoric systems versus general integrable
     systems}. The strategy to prove the inverse spectral conjecture
   for general integrable systems should be quite close to the
   semitoric case that we describe below.

   The only reason why one can do this for semitoric systems is
   because we have a complete symplectic theory in which a semitoric
   system is determined by some computable symplectic invariants.
   With our approach, having these invariants is a must, and for the
   moment we only have them for semitoric systems --- in the previous
   section we described a program to build symplectic invariants for
   more general systems.
 \item[(iii)] \emph{Local and semilocal inverse spectral
     formulas}. The challenge of proving a global result such as the
   inverse spectral conjecture passes by first proving local and
   semilocal inverse spectral results for singularities (focus-focus,
   elliptic, and hyperbolic), involving Bohr-Sommerfeld rules.  This
   can be avoided in the case of toric systems, where a global
   symplectic description of the system exists in terms of a single
   invariant: the image of the system, which is a polytope.
 \end{itemize}

 For each of the subsections \ref{sec:proper}, \ref{sec:higher},
 \ref{sec:antiperiodic}, and \ref{sec:degenerate} tentatively leading
 to a symplectic classification of integrable systems for the
 particular type of systems covered therein, one can try to push the
 strategy of proof outlined below to prove the inverse spectral
 conjecture for that type of systems.

 Unlike for Theorems \ref{inventiones} and \ref{acta}, we have not yet
 carried out the proof sketched below for semitoric systems, so the
 arguments are not as precise (and of course not guaranteed to
 work). However, we believe that this outline can develop into an
 efficient proof strategy.

\textup{\,}
\\
{\emph{Proof Strategy for Conjecture \ref{con}.}}
We describe the strategy for the case of semitoric systems only.  Even
though the setting to answer this question, in general, is the theory
Toeplitz operators, let's illustrate how to do this when $\hat{J},\,
\hat{H}$ are pseudodifferential operators and the phase space $(M,\,
\omega)$ is a cotangent bundle. The general strategy we present should
be adaptable to the general case of Toeplitz operators but this should
be technically more complicated and we will comment on this throughout
the proof.
\\
\\
{\bf \emph{Part A}} (Compute symplectic invariants from spectrum).  To
recover the number of singularities invariants, the volume invariant,
and the polygon invariant we equivalently recover: the image $F(M)$
together with the $J$\--coordinates of the focus\--focus values, and
the affine structure of $F(M)$.
\\
\\
\emph{Step 1} \emph{(Recovering $F(M)$ and $J$\--values of
  focus\--focus points)}.  We are able to do this because we have the
local density of states, i.e., the number of points of the spectrum in
a small ball of radius $\hbar^{\delta}$ multiplied by $\hbar^2$ and
divided by the volume of the ball (we take the ball of radius
$\hbar^{\delta}$, where $\delta \in (0,\,1)$; we are choosing
$\hbar^{\delta}$, but we could choose anything that goes to $0$
smoothly and slowly.)

The key point is that \emph{for pseudodifferential operators} we have
Bohr\--Sommerfeld rules and, as a consequence of these rules, we have
that, as $\hbar \to 0$, at a regular point the density function tends
to a smooth function of the center of a ball.  At an elliptic point,
the density function tends to a piecewise smooth but discontinuous
function. At a focus\--focus point, the density function tends to
infinity. To make this rigorous we have to work as in V\~ u Ng\d oc's
recent paper \cite{san-inverse}, and look at the local density of
states to find $F(M)$, the set of critical values, and the set of
non\--values. This method works, in general, for all integrable
systems, not necessarily semitoric.  For semitoric systems it may be
possible to even get an easier answer by applying the method line by
line for each line $J=\textup{constant}$ to cover the entire plane and
then analyze the eigenvalue spacing.

\begin{figure}[h]
  \begin{center}
    \includegraphics[width=0.4\linewidth]{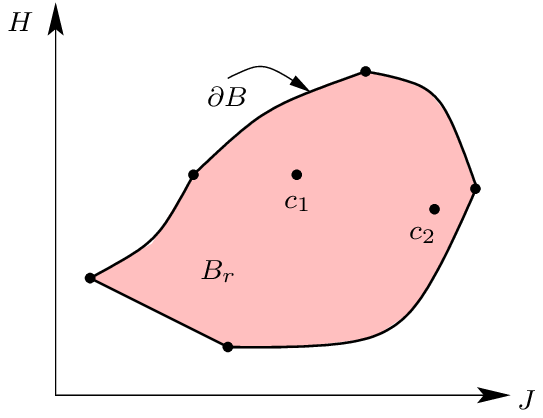}
    \caption{Image under $(J,\,H)$ of $M$ with focus\--focus values
      $c_1,\,c_2$.  $B=F(M)$ and $B_r$ is the set of regular values of
      $F$.}
    \label{fig:im2}
  \end{center}
\end{figure}

The singular Bohr-Sommerfeld rules would give the required result, in
case $\hat{J}$ and $\hat{H}$ were pseudodifferential operators. Of
course they are \emph{not}, since the phase space is not usually a
cotangent bundle. However $\hat{J},\, \hat{H}$ are semiclassical
Toeplitz operators and it is known that the algebra of Toeplitz
operators is ``microlocally equivalent'' to the algebra of
pseudodifferential operators. This works when $\hat{J},\, \hat{H}$ are
(self\--adjoint) Toeplitz operators.  What we have to do is to prove
the Bohr\--Sommerfeld rules for two degrees of freedom for general
Toeplitz operators. The one\--dimensional case has been done recently
by Charles (for Toeplitz operators; for pseudodifferential
operators it is well\--known). Because this is completely microlocal,
there is no expected difficulty in proving this (at least in
theory. But precise formulas are not easy to get; for instance, there
is no Maslov bundle in the Toeplitz case).
\\
\\
\emph{Step 2} \emph{(Recovering the affine structure of $F(M)$)}.  The
polygon invariant and the singular foliation type invariant (a Taylor
series on two variables) can both be defined in terms of the behavior
of the affine structure of $F(M)$ at its boundary. The affine
structure is itself defined by action integrals.  Therefore, it should
be possible to recover these invariants as soon as one can recover the
action integrals.

Because Bohr\--Sommerfeld rules give the spectrum in terms of action
integrals, they can be used, conversely, to compute actions from the
spectrum (see Figure \ref{fig:im}).

\begin{figure}[h]
  \begin{center}
    \includegraphics[width=0.5\textwidth]{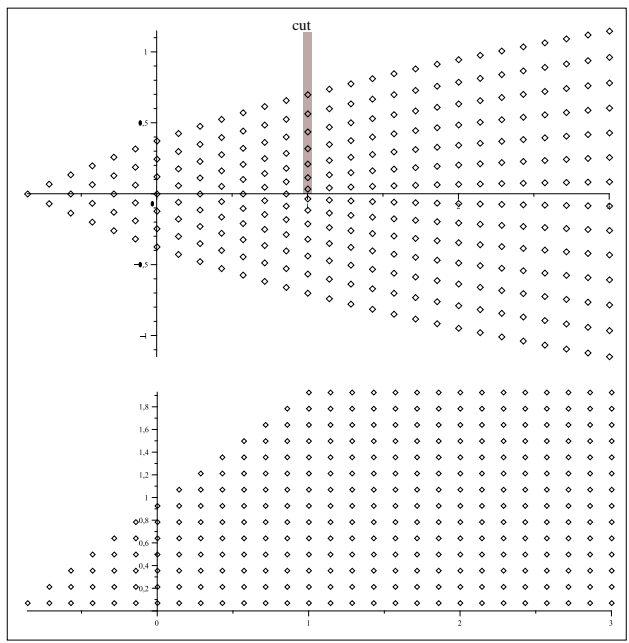}
    \caption{The joint spectrum of the quantum coupled
      spin\--oscillator has an integral affine structure on the
      regular values, which can be extended to the boundaries. Except
      along a vertical cut through the focus-focus critical value, one
      can develop this affine structure such that the joint
      eigenvalues become elements of $\hbar \Z^2$. The top picture is
      the joint spectrum of the quantum spherical pendulum.  At the
      bottom, we have developed the joint eigenvalues into a regular
      lattice. The number of eigenvalues in each vertical line is the
      same in both pictures.  The convex hull of the resulting set
      converges, as $\hbar \to 0$, to the semitoric polygon
      invariant.}
    \label{fig:im}
  \end{center}
\end{figure}

More precisely, the Bohr\--Sommerfeld rules should give the affine
structure in the regular part (=set of regular values). The singular
part of the affine structure is then obtained by taking limits of the
regular part. Let $(\gamma_1(c),\, \gamma_2(c))$ be a basis of
$\op{H}_1(F^{-1}(c),\, \mathbb{Z})$ depending smoothly on $c$. If $c_0
\in F(M)$ is a regular point, then there exists a ball $B$ around
$c_0$ such that
$$
\Big(\textup{JointSpec}(\hat{J},\, \hat{H}) \cap B \Big) =
f_{\hbar}(\hbar \mathbb{Z}^2)
$$
with
$$
f_{\hbar}(c)=f_0(c)+ \sum_{k\ge 1} \hbar^k f_k(c),
$$
where $f_0$ is a local diffeomorphism of $\mathbb{R}^2$,
\begin{eqnarray} \label{actionintegrals} f_0^{-1}(c):=\mathcal{A}(c)=
  \left(\mathcal{A}_1(c),\, \mathcal{A}_2(c)\right)
  =\left(\int_{\gamma_1(c)} \alpha,\,\, \int_{\gamma_2(c)}\alpha
  \right),
\end{eqnarray}
and $\alpha$ is any $1$\---form on some neighborhood of the fiber
$\Lambda_c$ containing $c$ such that $\op{d}\!\alpha=\omega$.  The
action integrals in (\ref{actionintegrals}) have a logarithmic
divergence at the focus\--focus singularities.

We claim that if one knows
$$
\textup{JointSpec}(\hat{J},\, \hat{H}) \cap F(M) \op{mod}(\hbar^2),
$$
then one should know $ f_{\hbar}(\hbar \mathbb{Z}^2) \op{mod}
\mathcal{O}(\hbar^2)$.  Of course, one needs to recover the whole
function, not just the discrete set, so this implies taking limits. A
crucial remark is that joint eigenvalues as $\hbar \to 0$ can approach
any arbitrary value in $F(M)$; from this we should be able to recover
$ f_{\hbar} \,\, \op{mod}\,\, \mathcal{O}(\hbar^2).  $ In particular,
the spectrum should completely determine $f_0$, and thus its inverse,
thereby giving the action integrals.

\emph{Step 3} (Recovering the singular foliation type
invariant). Let's try to recover the Taylor series $S^{\infty}$ at
focus\--focus singularities.  When we know the action\--integrals, in
order to recover the Taylor series, we use the formula:
$$
(S)^{\infty}=\textup{Taylor Series of the sum}
\underbrace{\mathcal{A}(c)}_{\textup{we know this}}\,\, \pm \,\,\,
\tilde{c} \op{ln}(\tilde{c}).
$$ 
We get $\mathcal{A}(c)=(\mathcal{A}_1(c),\,
\mathcal{A}_2(c))=\mathcal{A}_1(c)+\op{i} \mathcal{A}_2(c)$ from the
spectrum, which is the original data we are given.  The only
difficulty is that $\tilde{c}$ has to be expressed in Eliasson's
normal coordinates, so we need the \emph{complete} change of
coordinates.

We know that $\mathcal{A}(c)$ has a logarithmic behavior and it is
unique.  From the theory we know that there exists $\tilde{c}=g(c)$,
with $g$ a local diffeomorphism of $(\mathbb{R}^2,\, c_j)$ such that
$$
\underbrace{\overbrace{\mathcal{A}(c)}^{\mathcal{A}_1(c)+\op{i}\mathcal{A}_2(c)}\pm
  \overbrace{\tilde{c}\op{ln}(\tilde{c})}^{\tilde{c}(\op{ln}|\tilde{c}|
    +\op{i}\op{arg})}}_{\textup{We know}\,\, c,\,\, \textup{we don't
    know}\,\, \tilde{c} } \,\, \textup{is}\,\, \op{C}^{\infty} \,\,
\textup{on}\,\,\tilde{c}.
$$

Now, how many such $g$ are there? We claim that there exists a
\emph{unique} function $g$ (unique in the sense of its Taylor series
being unique) such that $\mathcal{A}(c)\pm
\tilde{c}\op{ln}(\tilde{c})$ is smooth.

The strategy to recover the last invariant, the twisting-index
invariant, should become more clear once the proof strategy above is
carried out. At the moment we have not made progress on this.
\\
\\
{\bf Part B} (Use the symplectic theory in \cite{PeVN2009,PeVN2011}
and Part A to recover integrable system). Once we have computed the
symplectic invariants, we can symplectically recover the integrable
system by Theorem \ref{inventiones} and Theorem \ref{acta}.

This concludes the comments on the general strategy to prove the
conjecture.

\subsection{\textcolor{black}{Conjecture \ref{con} in the case of
    toric systems}}

In the case of toric systems the only invariant is the polytope.  We
have recently proved \cite{ChPeVN2012} Conjecture \ref{con} for toric
integrable systems.

Recall from Section \ref{go} the geometric quantization procedure
introduced by Kostant and Souriau.  Let $(M,\, \omega)$ be a compact symplectic
manifold.  Given a \emph{prequantum} bundle $\mathcal{L} \rightarrow
M$, that is a Hermitian line bundle with curvature $\frac{1}{i} \om$
and a complex structure $j$ compatible with $\omega$, one defines the
quantum space as the space
$\mathcal{H}_k:=\mathrm{H}^0(M,\mathcal{L}^k)$ of holomorphic sections
of tensor powers $\mathcal{L}^k$ of $\mathcal{L}$, where $\hbar=1/k$
(when there is such a bundle, we say that $M$ is
\emph{prequantizable}).  Associated to such a quantization there is an
algebra $\mathscr{T} (M,{\mathcal{L}} , j)$ of Toeplitz operators
$(T_k \colon \Hilbert_k \rightarrow \Hilbert_k )_{k \in
  \mathbb{N}^*}$. Two Toeplitz operators $(T_k)_{k \in \N^*}$ and
$(S_k)_{k \in \N^*}$ \emph{commute} if $T_k$ and $S_k$ commute for
every $k$.

\begin{figure}[h]
  \centering
  \includegraphics[width=0.8\textwidth]{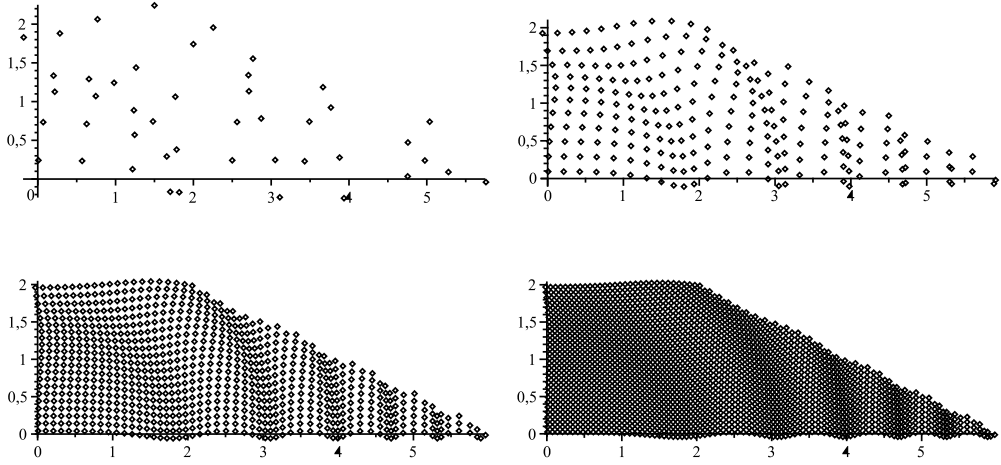}
  \caption{Sequence of images of the spectra of a quantum toric
    integrable system as the spectral parameter $\hbar$ goes to
    $0$. In the Hausdorff limit, corresponding to $\hbar=0$, the
    spectra converge to a polytope. The spectra lie in a plane, so
    they correspond to a four-dimensional integrable system with two
    degrees of freedom.}
\end{figure}

\begin{theorem}[Charles-Pelayo-V\~ u Ng\d oc]
  \label{theo:inverse-spectral}
  Let $(M, \, \omega, \, \mu : M \rightarrow \R^n)$ be a toric
  integrable system equipped with a prequantum bundle ${\mathcal{L}}$
  and a compatible complex structure $j$.  Let $T_1,\dots, T_n$ be
  commuting Toeplitz operators of $\mathscr{T} (M,{\mathcal{L}} , j)$
  whose principal symbols are the components of $\mu$. Then
$$ \De : = \lim_{k \rightarrow \infty}  \textup{JointSpec}(T_1,\dots, T_n) $$
is a Delzant polytope and $(M, \om, \mu)$ is isomorphic with
$(M_{\Delta}, \om_{\Delta}, \mu_{\Delta})$.

In other words, one can recover the classical system from the limit of
the joint spectrum.
\end{theorem}

The proof of this theorem combines geometric and algebraic methods,
following the framework introduced in Duistermaat\--Pelayo
\cite{DuPe2009}, with semiclassical analytic techniques inspired by
recent work of Charles \cite{Ch2003a,  Ch2006b}.

\subsection{\textcolor{black}{Existence of quantization}}
The quantum system is the data for us; right now we are not interested
in the question of whether one can find, or how to find, such
operators -- that is a different question. However, in the case of
toric systems, we could address this question, as a consequence of a
global normal form theorem we proved.
 
Our conclusion was that if $(M, \, \omega)$ is a prequantizable
compact connected symplectic manifold of dimension $2n$ equipped with
a toric integrable system given by a momentum map $\mu:=(\mu_1,\,
\ldots,\,\mu_n) \colon M \to \mathbb{R}^n$ of a Hamiltonian
$n$\--torus action on $M$, then for any Toeplitz quantization
$(\mathcal{H}_k)_{k\in \NM^*}$ of $M$ there exists a quantization of
$\mu$, i.e., a set of $n$ Toeplitz commuting operators whose principal
symbols are $\mu_1,\dots,\mu_n$. In the general analytic case, the
obstruction has been found by Garay.

\section{Other projects}

\subsection{Algorithmic version of spectral goal and
    proof strategy}
\label{sec:algo}

Given the semiclassical joint spectrum of a quantum system, we would
like to have an algorithm that implements the strategy of proof of the
inverse spectral conjecture to first compute the symplectic invariants
from the spectrum (Step 1) and then constructs the system from the
invariants (Step 2). In the case of toric systems (for which the
inverse spectral conjecture is a theorem), the only invariant is (by
Delzant's theorem) the polytope and the recipe to recover this
polytope from the spectrum is already described in Theorem
\ref{theo:inverse-spectral}: take the Hausdorff limit of the sequence
of spectra; this gives Step 1.

It is not completely clear how to implement an algorithm to carry out
Step 2, which is Delzant's construction (an explicit construction by
symplectic reduction on some complex space $\C^N$). From the point of
view of applications, it may be enough to extract information from the
polytope about $M$ and give a list of its properties (number of fixed
points, symplectic volume, etc.) rather than try to give a
construction of the toric system, which may be of little use. The key
point is that \emph{all} the information about the toric system is in
the polytope, in the same way that all the information about a
semitoric system is in the list of our five symplectic invariants.

\subsection{\textcolor{black}{Type theory and Univalent Foundations}}

The formalization of mathematics in type theory, and in particular in
the Coq system, in the language and setting of Voevodsky's Univalent
Foundations \cite{[1]}, and using Voevodsky's libraries is a project
that one could view as a modern version of Bourbaki.  The Univalent
Foundations are constructive and new mathematics is likely to emerge
from the necessity to formulate various components of the integrable
systems theory constructively.

The goal is to formalize the theory of finite dimensional integrable
systems.  We are suggesting that this is a promising way to think
about the organization of the theory of integrable systems, which
should facilitate progress in the field in the future. In particular,
one expects to create novel mathematics by means of this
formalization.  In short, the idea of the Voevodsky's program is to
develop mathematics within the world of homotopy types.  Voevodsky's
program uses deep insights and tools from homotopy theory. The project
of \'A. Pelayo, M.A. Warren, and V. Voevodsky has the following steps:

\begin{itemize}
\item[i)] \emph{formalize some basic $p$\--adic analysis} (a
  preliminary file dealing with this case may be found in \cite{[2]});

\item[ii)] \emph{to formalize $p$\--adic integrable systems}. In fact
  first one needs to rigorously give the ``correct" formulation of
  $p$\--adic integrable systems.  This is a very new and interesting
  topic on its own (separately from type theory and Coq), as a problem
  about integrable systems. In fact, it is non\--trivial to describe a
  useful notion of $p$\--adic integrable system.

\item[iii)] \emph{continue to formalize general integrable systems}.
\end{itemize}

We give three reasons for which we believe one should be interested in
this formalization:

\begin{itemize}
\item[a)] \emph{Creation of Algorithms}: because everything is done in
  a system that has good algorithmic properties (assuming that
  Voevodsky's conjecture regarding the constructive status of the
  Univalence Axiom holds), {one should be able to extract good
    algorithms from the proofs}. This can be extraordinarily useful,
  because we may potentially be able to get explicit information when
  otherwise we would only have existence results;

\item[b)] \emph{Numerical information about experiments (``Applied
    Spectral Goal")}: this section should have applications to
  Section~\ref{sec:algo}.  The construction of algorithms in a),
  should help with the following outstanding problem, of which we
  discussed an exact version earlier in this paper: given numerical
  spectral data about a quantum integrable system (coming from an
  experiment), extract (optimizing the algorithms obtained as proof
  terms in (a)) an algorithm to reconstruct the classical integrable
  system.

\item[c)] \emph{New theoretical insights}: by virtue of the fact that
  the project involves formalizing mathematics in the world of
  homotopy, the notions formalized can be seen in many cases as ``up
  to homotopy versions" of the original notions, and therefore are
  more general.  As for example derived algebraic geometry suggests, it
  is often mathematically more useful to work with mathematical
  structures up to homotopy than it is to work with their strict
  analogues. The hope is therefore that this project will result in
  the creation of novel mathematics and not just the formalization of
  known mathematics.
\end{itemize}

\section{Perspective of this paper and closing
    remarks} Not only the field of integrable systems is vast on its own right, but
also it is rare to find an area of mathematics in which they do not
occupy a privileged position: from symplectic and algebraic geometry
to representation theory, microlocal analysis, spectral theory,
partial differential equations, numerical analysis, etc.

\subsection{\textcolor{black}{Symplectic theory: perspective and
    connections}} \label{sec:r}

Although this paper is not a survey, to help the interested gain some
perspective, we have included in the references a small selection of
the important contributions to integrable Hamiltonian systems which
are most directly relevant to this paper. We also try to briefly
outline a few connections to other works, without claiming any
completeness.

\subsubsection{\textcolor{black}{Glimpse of references}}
 
For the reader's aide we have cited:
\begin{itemize}
\item[(i)] the essential works on which the ideas proposed in this
  paper are based;
\item[(ii)] works in progress on open problems outlined in the paper;
\item[(iii)] a few works that give perspective on the current paper
  and the subjects on which it touches.
\end{itemize}
See our articles \cite{PeVN2012, ChPeVN2012}, and the references
therein, for further information on related works and other
motivations.  The first few sections of \cite{PeVN2012} provide a
brief review of several fundamentals results of the theory of finite
dimensional integrable systems due to Arnold, Atiyah,
Carath{\'e}odory, Darboux, Delzant, Duistermaat, Dufour, Eliasson,
Guillemin, Liouville, Mineur, Molino, Sternberg, Toulet and
N. T. Zung, among others, which have been key ingredients in the most
recent developments outlined in the present paper.

Finally, we would like to point out that in the particular case that
$M$ is a cotangent bundle, fundamental work on integrable systems was
done by Hitchin \cite{Hi1987}. The work of Lax \cite{lax} has also
been extremely influential.

\subsubsection{\textcolor{black}{Atiyah, Delzant, Guillemin, Kirwan,
    and Sternberg theory}}

As indicated in Section \ref{sec:setting}, the authors' intuition on
integrable systems has been guided by a remarkable theory mostly
developed in the 70s and 80s by Kostant, Atiyah \cite{atiyah},
Guillemin\--Sternberg \cite{gs}, Delzant \cite{delzant} and Kirwan
\cite{kirwan} in the context of Hamiltonian torus actions.

The construction of the convex polytope $\Delta$ in the Atiyah and
Guillemin\--Sternberg Theorem (Theorem \ref{theo:ags}) was the
motivation and driving force behind the construction of the ``polygon
invariant'': item (3) in Theorem \ref{inventiones}.  Indeed,
Hamiltonian $n$\--torus actions on symplectic $2n$\--manifolds form an
important class of completely integrable systems, with well\--behaved
singularities, usually referred to as \emph{toric systems}.

Delzant built on Theorem \ref{theo:ags} to give a classification of
toric systems. His theorem, which in dimension $4$ is a particular
case of Theorems \ref{inventiones}, \ref{acta}, was a main motivation
for these results.

\subsubsection{\textcolor{black}{Kolmogorov\--Arnold\--Moser theory}}

One motivation to study integrable systems comes from
Kolmogorov\--Arnold\--Moser (KAM) theory.  Since integrable systems
are ``solvable'' in a precise sense, one expects to find valuable
information about the behavior of dynamical systems that are obtained
by \emph{small perturbations} of them, and then the powerful KAM
theory comes into play to deal with the properties
of the perturbations (persistence of quasi\--periodic motions).

An important limitation of KAM techniques is that they require the
unperturbed integrable system to be written in action-angle
coordinates.  Having global action-angle coordinates is very
exceptional; thus, most applications of KAM are limited to
neighborhoods of regular Lagrangian tori.

Meanwhile, it has become clear that singularities of integrable
systems may enter KAM theory~: see for
instance~\cite{zung-kolmogorov} where the various
Kolmogorov conditions are deduced from the existence of hyperbolic or
focus-focus singularities.  This opens the way to a global KAM
theory. Recent results of~\cite{broer-al} show how to construct Cantor
sets of invariant tori near focus-focus singularities, and thus define
Hamiltonian monodromy for the perturbed system. It is tempting to
apply these ideas to semitoric systems, and see which parts of the
integral affine structure (and thus, which symplectic invariants of
Section~\ref{sec:semitoric}) survive the perturbation.

If this can be achieved, the quantum counterpart becomes fascinating:
the reminiscence of the integral affine structure should have an
effect on the spectrum of the perturbed quantum Hamiltonian. What
makes it mysterious is that one doesn't have any joint spectrum in the
sense of Section~\ref{sec:spectral} for the perturbed
Hamiltonian. This issue might perhaps be overcome by switching to
non-selfadjoint perturbations as in~\cite{san-hitrik-sjoestrand}.

\subsubsection{\textcolor{black}{Fomenko School theory}}

In the 1980s and 1990s, the Fomenko school developed a far reaching
Morse theory for regular energy surfaces of integrable systems, which
is related, and serves as an inspiration, for several of the problems
(for instance Problems \ref{pr1a}, \ref{pr1b}, \ref{pr1b2},
\ref{pr1c}) discussed in Section \ref{sec:topology}.
 
One of our motivations to study integrable systems comes from the
theory of singularities of fibrations $$\xi \colon M \to
\mathbb{R}^n$$ that they developed \cite{bolsinov-fomenko-book}.  They
gave an extensive treatment of the topological properties of
integrable systems viewed as fibrations over $\mathbb{R}^n$, and their
work has exerted a decisive influence on the problems presented in
Section \ref{sec:topology} of the present paper.

\subsubsection{\textcolor{black}{Singular affine structures}}

A semitoric system as in Definition
\ref{def:semitoric} gives rise to a torus fibration with
singularities, and its base space becomes endowed with a singular
integral affine structure.
Singular affine structures are of key importance in various parts of
symplectic topology, mirror symmetry, and algebraic geometry -- for
example they play a central role in the work of 
Gross and Siebert \cite{gross, grs3}, Kontsevich and Soibelman
\cite{KS}, and Symington \cite{s}.  

In the semitoric case, the singular affine structure is encoded in the
first three invariants/ingredients in Theorems \ref{inventiones},
\ref{acta}:
\begin{itemize}
\item[{\rm (1)}] \emph{number of singularities}: the number of
  focus\--focus singularities $m_f$;
\item[{\rm (2)}] \emph{singular foliation type}: a formal Taylor
  series $S(X,Y)$ at each focus\--focus singularity;
\item[{\rm (3)}] \emph{polygon invariant}: a class of polygons
  equipped with $m_f$ oriented vertical lines (see Figure
  \ref{fig:weightedpolygon});
\end{itemize}

\subsubsection{\textcolor{black}{Quasistates}}

Theorems \ref{inventiones} and \ref{acta} characterize semitoric
systems and show that the collection of such systems is huge.  Simple
examples of semitoric semitoric systems appear in the theory of
symplectic quasi-states, see the article Eliashberg\--Polterovich
\cite[page 3]{eliashberg}.  Semitoric systems could provide many other
interesting examples for the study of displaceability and related
questions in that context.

\subsubsection{\textcolor{black}{Mechanics}}

As described in Sections \ref{sec:intro} and \ref{sec:symplectic}, an
important motivation for the present paper has been to provide a
unified framework to study integrable systems. Many integrable systems
of fundamental importance arise in classical mechanics: for instance,
the \emph{coupled spin\--oscillator} (also called the
\emph{Jaynes\--Cummings model}), the \emph{spherical pendulum}, the
\emph {two\--body problem}, the \emph{Lagrange top}, the \emph{three
  wave interaction} (see Section \ref{sec:first} for more details
about these systems). In fact, more generally, integrable systems can
also be found in the theory of geometric phases, rigid body systems,
elasticity theory and plasma physics, and have been extensively
studied by many authors.

Integrable systems appear in numerous contexts, which we have not
explored in the present paper. N. Reshetikhin's survey lectures on
aspects of classical and quantum integrability \cite{Re2010} are an
excellent source for related important developments; the lectures,
among other contributions, outline the relation between solvable
models in statistical mechanics, and classical and quantum integrable
spin chains, which are closely related to the spin-oscillator example
in Section \ref{sec:first}.

\subsection{\textcolor{black}{Spectral theory: perspective and
    connections}}

This part of the program outlined in this paper combines geometric
ideas in the complex and symplectic settings with microlocal analytic
methods dealing with semi-classical pseudodifferential and Toeplitz
operators.

\subsubsection{\textcolor{black}{Isospectrality for integrable
    systems}}

The central theme surrounding the ``Spectral Goal for Integrable
Systems" (Sections~\ref{sec:intro} and \ref{sec:spectral}) is the
\emph{isospectrality} question for quantum integrable systems. In
other words, in any finite dimension: does the semiclassical joint
spectrum of a quantum toric integrable system, given by a sequence of
commuting Toeplitz operators acting on quantum Hilbert spaces,
determines the classical system given by the symplectic manifold and
Poisson commuting functions, up to symplectic isomorphisms? This type
of symplectic isospectral problem belongs to the realm of classical
questions in inverse spectral theory and microlocal analysis, going
back to pioneer works of Colin de Verdi{\`e}re on cotangent bundles
\cite{colinII, CdV2} and Guillemin\--Sternberg \cite{GuSt} in the
1970s and 1980s.

\subsubsection{\textcolor{black}{Microlocal analysis of integrable
    systems}}

The notion of a quantum integrable system
(Definition~\ref{defi:quantum}), as a maximal set of commuting quantum
observables, dates back to the early quantum mechanics, to the works
of Bohr, Sommerfeld and Einstein~\cite{E1917}.

However, the most basic results in the symplectic theory of classical
integrable systems like Darboux's theorem or action-angle variables
could not be used in Schr{\"o}dinger's quantum setting at that time
because they make use of the analysis of differential (or
pseudodifferential) operators in phase space, known now as microlocal
analysis, which was developed only in the 1960s.

The microlocal analysis of action-angle variables starts with the
works of Duistermaat~\cite{Du1974} and Colin de Verdi{\`e}re
\cite{colinII,CdV2}, followed by the semiclassical theory by
Charbonnel \cite{[5]}, and more recently by V\~{u} Ng\d{o}c
\cite{san-spectral},  Toth and Zelditch \cite{[50]}, Charbonnel
and Popov \cite{ChPo1999}, Melin\--Sj{\"o}strand \cite{MeSj2003}, and
many others.
 
Effective models in quantum mechanics often require a compact phase
space, and thus cannot be treated using pseudodifferential
calculus. For instance the natural classical limit of a quantum spin
is a symplectic sphere.  The study of quantum action-angle variables
in the case of compact symplectic manifolds treated in this paper was
started by Charles \cite{Ch2003a}, using the theory of Toeplitz
operators.

\subsubsection{\textcolor{black}{Kac's question}}

The \emph{Spectral Goal} (Sections~\ref{sec:intro} and
\ref{sec:spectral}) fits in the framework of ``isospectral questions":
what is the relation between two operators that have the same
spectrum? The question of isospectrality has been considered by many
authors in different contexts, and may be traced back to a more
general question of H. Weyl.  In the case of the
Riemannian Laplacian, the question is
perhaps most famous thanks to Kac's article \cite{Ka66} (who
attributes the question to S. Bochner), which also popularized the
phrase: ``can one hear the shape of a drum?".

Bochner and Kac's question has a negative answer in this case, even
for planar domains with Dirichlet boundary conditions (which is the
original version posed in \cite{Ka66}). As P. Sarnak mentioned to us,
a much better question to ask is whether the set of isospectral
domains is finite. There are many related works, see for instance
Milnor \cite{Mi1964},  Osgood\--Phillips\--Sarnak
\cite{OsPhSa1989}, Buser \cite{Bu86} and Gordon\--Webb\--Wolpert \cite{GoWeWo92}.

Inverse type results in the realm of spectral geometry have been
obtained by many other authors, see for instance Br{\"u}ning\--Heintze
\cite{Br1984b}, Colin de Verdi{\`e}re \cite{colinII, CdV2}, and Zelditch \cite{Ze2009},
and the references therein.

\subsubsection{\textcolor{black}{Quantum spectroscopy}}

>From the point of view of applications, one of the main goals of the
program described in the present paper is to solve concrete problems
arising in quantum molecular spectroscopy (see Section~\ref{sec:algo}). 
Indeed, physicists and chemists were the first to become interested in
semitoric systems which appear naturally in the context of quantum
chemistry. Many groups have been working on this topic, to name a few:
M. Child's group in Oxford, J. Tennyson's at University College
London, F. De Lucia's at Ohio State University, B. Zhilinski\'{\i}'s
at Dunkerque, and M. Joyeux's at Grenoble. See \cite{Fi2009, Sa1995,
  Ch1999}.  The question which arose from the work of these physicists
and chemists is whether one can give a finite collection of invariants
characterizing systems of this nature (Section~\ref{sec:semitoric}).
We are very grateful to B. Zhilinski\'i. for bringing this question to
our attention.

\subsection{\textcolor{black}{Conclusion}}
In this paper we have suggested some ideas to work towards developing
the symplectic and spectral theory of finite dimensional integrable
Hamiltonian systems. The paper has laid out many open problems.  

This paper has the authors' articles \cite{PeVN2009,PeVN2011} as
reference points, and it describes several somewhat disjoint paths in
which one can make contributions to the topics at hand by dealing with
the difficulties one at a time.  We expect that each of these disjoint
paths will involve substantial challenges. However, we believe that
the overall strategy stands a reasonable chance of success, by what we
mean that it should lead to improving our knowledge of the symplectic
and spectral aspects of integrable systems significantly.

We expect that further information about the program described in this
paper will be available in the authors' websites listed in this paper.

\section*{Acknowledgments} 
We are indebted to Tudor Ratiu for his extraordinary dedication to
reading this paper.  His numerous comments prompted us to rethink
parts of the paper and have greatly improved the exposition and the
clarity.  The authors are grateful to Helmut Hofer for his essential
support that made it possible for VNS to visit the Institute for
Advanced Study several times in 2011.  We are also thankful to the
referees for their careful reading and many suggestions, which have
improved the paper significantly.  We are grateful to Jochen
Br\"uning, Nigel Hitchin, Helmut Hofer, Darryl Holm, Tudor Ratiu,
Peter Sarnak, Thomas Spencer, Vladimir Voevodsky, and Michael
A. Warren for helpful discussions.  We are grateful to Alan Weinstein
for comments on a preliminary version.  AP was partly supported by an
NSF Postdoctoral Fellowship, NSF Grants DMS-0965738 and DMS-0635607,
an NSF CAREER Award, an Oberwolfach Leibniz Fellowship, Spanish
Ministry of Science Grant MTM 2010-21186-C02-01, and by the Spanish
National Research Council. VNS was supported by the Institut
Universitaire de France and a NONAa grant from the French ANR.

\vskip3em
\noindent
\\
{\bf {\'A}lvaro Pelayo} \\
School of Mathematics\\
Institute for Advanced Study\\
Einstein Drive\\
Princeton, NJ 08540 USA.
\\
\\
\noindent
Washington University,  Mathematics Department \\
One Brookings Drive, Campus Box 1146\\
St Louis, MO 63130-4899, USA.\\
{\em E\--mail}: \texttt{apelayo@math.wustl.edu} \\
{\em Website}: \url{http://www.math.wustl.edu/~apelayo/}\\

\medskip\noindent

\noindent
\noindent
{\bf San V\~u Ng\d oc} \\
Institut Universitaire de France
\\
\\
Institut de Recherches Math\'ematiques de Rennes\\
Universit\'e de Rennes 1\\
Campus de Beaulieu\\
F-35042 Rennes cedex, France\\
{\em E-mail:} \texttt{san.vu-ngoc@univ-rennes1.fr}\\
{\em Website}: \url{http://blogperso.univ-rennes1.fr/san.vu-ngoc/}

\end{document}